\documentclass[MSNbibl,number,citesort,dvips]{arxbj}
\usepackage{upgreek}

\aid{0}
\volume{19}
\issue{5B}
\pubyear{2013}
\firstpage{2494}
\lastpage{2523}
\doi{10.3150/12-BEJ460} 

\makeatletter
\newcommand{\ppi}{\uppi}
\newcommand{\rrvert}{\vert}
\newcommand{\llvert}{\vert}
\newcommand{\overset}{\stackrel}
\newtheorem{theorem}{Theorem}
\newtheorem{corollary}[theorem]{Corollary}
\newproclaim{definition}[theorem]{Definition}
\newtheorem{lemma}[theorem]{Lemma}
\newtheorem{proposition}[theorem]{Proposition}
\newremark{remark}[theorem]{Remark}
\newcommand{\sign}{\operatorname{sgn}}
\makeatother

\begin{document}
\begin{frontmatter}

\title{The Lamperti representation of real-valued self-similar Markov processes}
\runtitle{The Lamperti representation of real-valued self-similar
Markov processes}

\begin{aug}
\author[1]{\fnms{Lo\"{i}c} \snm{Chaumont}\thanksref{1}\ead[label=e1]{loic.chaumont@univ-angers.fr}},
\author[2]{\fnms{Henry} \snm{Pant\'i}\corref{}\thanksref{2,e2}\ead[label=e2,mark]{henry@cimat.mx}} \and
\author[2]{\fnms{V\'ictor} \snm{Rivero}\thanksref{2,e3}\ead[label=e3,mark]{rivero@cimat.mx}}
\runauthor{L. Chaumont, H. Pant\'i and V. Rivero} 
\address[1]{LAREMA, D\'epartement de Math\'ematiques, Universit\'e
d'Angers. 2, Bd Lavoisier - 49045, Angers Cedex~01, France. \printead{e1}}
\address[2]{Centro de Investigaci\'on en Matem\'aticas (CIMAT A.C.),
Calle Jalisco s/n, 36240 Guanajuato, Guanajuato, M\'exico. \printead{e2,e3}}
\end{aug}

\received{\smonth{11} \syear{2011}}
\revised{\smonth{6} \syear{2012}}

%
\begin{abstract}
In this paper, we obtain a Lamperti type representation for real-valued
self-similar Markov processes, killed at their hitting time of zero.
Namely, we represent real-valued self-similar Markov processes as time
changed multiplicative invariant processes. Doing so, we complete Kiu's
work [\textit{Stochastic Process. Appl.} \textbf{10} (1980) 183--191],
following some ideas in Chybiryakov [\textit{Stochastic Process. Appl.} \textbf{116} (2006) 857--872] in
order to characterize the underlying processes in this representation.
We provide some examples where the characteristics of the underlying
processes can be computed explicitly.
\end{abstract}

%
\begin{keyword}
\kwd{Lamperti representation}
\kwd{L\'evy processes}
\kwd{multiplicative invariant processes}
\kwd{self-similar Markov processes}
\end{keyword}

\end{frontmatter}

\section{Introduction}

Semi-stable processes were introduced by Lamperti in \cite{Lamperti62}
as those processes satisfying a scaling property. Nowadays this kind of
processes are known as self-similar processes. Formally, a c\`adl\`ag
stochastic process $X=(X_t, t\geq0)$, with $X_0=0$, and Euclidean
state space $E$, is self-similar of order $\alpha>0$, if for every
$a>0$, the processes $(X_{at}, t\geq0)$ and $(a^{\alpha}X_t, t\geq
0)$, have the same law. Lamperti proved that the class of self-similar
processes is formed by those stochastic processes that can be obtained
as the weak limit of sequences of stochastic processes that have been
subject to an infinite sequence of dilations of scale of time and
space. More formally, the main result of Lamperti in \cite{Lamperti62}
can be stated as follows: let $(\widetilde{X}_t, t\geq0)$ be a
stochastic process defined in some probability space $(\Omega, \mathcal
{F}, \mathbb{P})$ with values in $E$. Assume that there exists a
positive real function $f(\eta) \nearrow\infty$ such that the process
$(\widetilde{X}_t^\eta, t\geq0)$ defined by
\[
\widetilde{X}_t^\eta= \frac{\widetilde{X}_{\eta t}}{f(\eta)}, \qquad t\geq0,
\]
converges to a non-degenerated process $X$ in the sense of
finite-dimensional distributions. Then, $X$ is a self-similar process
of order $\alpha$ and $f(\eta)=\eta^\alpha L(\eta)$, for some $\alpha
>0$, where $L$ is a slowly varying function. The converse is also true,
every self-similar process can be obtained in such a way.

If $X$ is a Markov process with stationary transition function
$P_t(x,A)$, then the self-similarity property written in terms of its
transition function takes the form
%
\begin{equation}
\label{eq111} P_{at}(x, A) = P_t \bigl(a^{-1/\alpha}
x, a^{-1/\alpha} A \bigr)
\end{equation}
for all $a>0$, $t\geq0$, $x\in E$, and all measurable sets $A$. We
will assume that $X$ is a strong Markov process and refer to it as a
\textit{self-similar Markov process of index} $\alpha>0$.

From now on, $\Omega$ denotes the space of c\`adl\`ag paths, $X$ the
coordinates process and $(\mathcal{F}_t, t\geq0)$ its natural
filtration, that is, $\mathcal{F}_t = \sigma(X_s, s\leq t)$.

There are many other ways than (\ref{eq111}) to define self-similar
Markov processes. The definition used in this paper is the following.

\begin{definition} \label{def1}
Let $E$ be $[0,\infty)$ or $\mathbb{R}^n$. We will say that $\{
X^{(x)}=(X, \mathbb{P}_x), {x\in E}\}$ is a family of $E$-valued
self-similar Markov processes with index $\alpha>0$ if it is a c\`adl\`
ag strong Markov family with state space $E$, and that satisfies that
for every $c>0$,
\[
\bigl\{(cX_{c^{-\alpha}t}, t\geq0), \mathbb{P}_x \bigr\} \overset
{\mathcal{L}}= \bigl\{(X_t, t\geq0), \mathbb{P}_{cx} \bigr
\} \qquad \forall x\in E.
\]
\end{definition}

The case $E=[0,\infty)$ was first investigated by Lamperti in \cite{Lamperti72} and has further been the object of many studies, see, for
instance, Bertoin  and Yor \cite{Bertoin-Yor}, Carmona, Petit and Yor \cite{Carmona-Petit-Yor} and the
reference therein. Here, we summarize some of his main results. Let $T$
be the first hitting time of zero for $X$, that is,
\[
T=\inf\{t>0\dvt X_t=0\},
\]
with $\inf\{\varnothing\}=\infty$. Then, for any starting point $x>0$,
one and only one of the following cases holds:
\begin{enumerate}[{C.1}]
\item[{C.1}] $T=\infty$, a.s.

\item[{C.2}] $T<\infty$, $X_{T-}=0$, a.s.

\item[{C.3}] $T<\infty$, $X_{T-}>0$, a.s.
\end{enumerate}
We refer to {C.1} as the class of processes that never reach
zero, processes in the class {C.2} hit zero continuously, and
those in the class {C.3} reach zero by a jump. In particular, if
$T$ is finite, then the process reaches zero continuously or by a jump.
Another important result in Lamperti \cite{Lamperti72} is the representation of
positive self-similar Markov processes as the exponential of L\'evy
processes time changed by the inverse of their exponential functional.
This representation is known as the Lamperti representation and its
extension to real-valued processes is one of the main motivations of
this paper. Formally, the Lamperti representation can be stated as
follows. Assume that the process $X$ is absorbed at $0$. Let $(\xi_t,
t\geq0)$ be the process defined by
\[
\exp\{\xi_t\} = x^{-1}X_{\nu(t)},\qquad  t\geq0,
\]
where
\[
\nu(t) = \inf \biggl\{s>0\dvt \int_0^s
(X_u)^{-\alpha} \,\mathrm{d}u > t \biggr\},
\]
with the usual convention $\inf\{\varnothing\}=+\infty$. Then, under
$\mathbb{P}_x$, $\xi$ is a L\'evy process. Furthermore, $\xi$ satisfies
either (i) $\limsup_{t\rightarrow\infty} \xi_t = \infty$ a.s.,
(ii) $\lim_{t\rightarrow\infty} \xi_t = -\infty$ a.s. or (iii) $\xi
$ is a L\'evy process killed at an independent exponential time $\zeta
<\infty$ a.s., depending on whether $X$ is in the class~{C.1},
{C.2} or {C.3}, respectively. Note that since an
exponential random variable with parameter $q$ is infinite if only if
$q=0$, then we can always consider the process $\xi$ as a L\'evy
process killed at an independent exponential time $\zeta$ with
parameter $q\geq0$. Conversely, let $(\xi, \mathbf{P})$ be a L\'evy
process killed at an exponential random time $\zeta$ with parameter
$q\geq0$ and cemetery point $\{-\infty\}$. Let $\alpha>0$ and for
$x>0$, define the process $X^{(x)}$ by
\[
X_t^{(x)} = x\exp\{ \xi_{\tau(t x^{-\alpha})} \}, \qquad t\geq0,
\]
where
\[
\tau(t) = \inf \biggl\{u>0\dvt \int_0^u \exp\{
\alpha\xi_s \}\,\mathrm{d}s > t \biggr\}.
\]
Then, $(X^{(x)})_{x>0}$ is a positive self-similar Markov process of
index $\alpha>0$ which is absorbed at~0. Furthermore, the latter
classification depending on the asymptotic behaviour of $\xi$ holds. An
important relation between $T$ and the exponential functional of the L\'
evy process $\xi$ is $(T, \mathbb{P}_x) \overset{\mathcal{L}}=
(x^\alpha\int_0^\zeta\exp\{\alpha\xi_s\} \,\mathrm{d}s, \mathbf{P})$. Further
details on this topic can be found in Lamperti \cite{Lamperti72}, Bertoin and Yor \cite{Bertoin-Yor}.

In Kiu \cite{Kiu80}, the case of $\mathbb{R}^n$-valued self-similar Markov
processes was studied. The main result in Kiu \cite{Kiu80} asserts that,
if $X$ killed at $T$ is a Feller self-similar Markov process, then the
process~$Y$ defined by
\[
Y_t = X_{\nu(t)},\qquad  t\geq0,
\]
where
\[
\nu(t) = \inf \biggl\{s>0\dvt \int_0^s
|X_u|^{-\alpha}  \,\mathrm{d}u > t \biggr\},
\]
is a Feller multiplicative invariant process, that is, $Y$ is a Feller
process with semigroup $Q_t$ satisfying
%
\begin{equation}
\label{eq222} Q_t(x, A) = Q_{t}(ax,aA)
\end{equation}
for all $x\neq0$, $a,t$ positive and $A \in\mathcal{B}(\mathbb
{R}^{n}\setminus\{0\})$. Another way to write (\ref{eq222}) is
\[
Q_t \bigl(x, a^{-1}A \bigr) = Q_{t} \bigl(|a|x,
\sign(a)A \bigr)
\]
for all $t$ positive, $x, a\neq0$ and $A \in\mathcal{B}(\mathbb
{R}^n\setminus\{0\})$. This property may also be written in terms of
the process $Y$ as follows:
%
\begin{equation}
\label{eq333} \bigl\{(aY_t, t\geq0), \mathbb{P}_{x}
\bigr\} \overset{\mathcal{L}}= \bigl\{ \bigl(\sign(a)Y_t, t\geq0
\bigr), \mathbb{P}_{|a|x} \bigr\}
\end{equation}
for all $x, a\neq0$. In Kiu \cite{Kiu80}, the converse of this result has
not been proved but using (\ref{eq333}), it is easy to verify that it
actually holds. Formally, let $Y$ be a strong Markov process taking
values in $\mathbb{R}^n\setminus\{0\}$ and satisfying (\ref{eq333}).
Let $\alpha>0$ and define the process $X$ by
\[
X_t = Y_{\varphi(t)},\qquad  t\geq0,
\]
where
\[
\varphi(t) = \inf \biggl\{s>0\dvt \int_0^s
|Y_u|^\alpha \,\mathrm{d}u>t \biggr\},
\]
with $\inf\{\varnothing\} = \infty$. Then $X$ is a $\mathbb{R}^n$-valued
self-similar Markov process of index $\alpha>0$ which is killed at $T$.
It is important to mention that no explicit form of $Y$ has been given
in Kiu \cite{Kiu80}. Giving a construction of Feller multiplicative
invariant processes taking values in $\mathbb{R}^\ast:=\mathbb
{R}\setminus\{0\}$, that we will call Lamperti--Kiu processes, is
another main motivation of this paper.

\begin{definition}
Let $Y=(Y_t, t\geq0)$ be a c\`adl\`ag process. We say that $Y$ is a
Lamperti--Kiu process if it takes values in $\mathbb{R}^\ast$, has the
Feller property and (\ref{eq333}) is satisfied.
\end{definition}

A subclass of Lamperti--Kiu processes has been studied by Chybiryakov in \cite{Chybiryakov}
who gave the following definition. Let $Y$ be a
$\mathbb{R}^\ast$-valued c\`adl\`ag process defined on some probability
space $(\Omega, \mathcal{F}, \mathbf{P})$ such that $Y_0 = 1$. It is
said that $Y$ is a multiplicative L\'evy process if for any $s,t>0$,
$Y_t^{-1}Y_{t+s}$ is independent of $\mathcal{G}_t = \sigma(Y_u, u\leq
t)$ and the law of $Y_t^{-1}Y_{t+s}$ does not depend on $t$. It can be
shown that if $Y$ is a multiplicative L\'evy process, then $Y$ is
Markovian and its semigroup satisfies (\ref{eq222}). Furthermore, there
exist a L\'evy process $\xi$, a Poisson process $N$ and a sequence
$U=(U_k, k\geq0)$ of i.i.d. random variables, all independent, such that
%
\begin{equation}
\label{eq444} Y_t = \exp \Biggl\{ \xi_t + \sum
_{k=1}^{N_t} U_k + \mathrm{i}\ppi N_t
\Biggr\}, \qquad t\geq0.
\end{equation}
The converse is also true, that is, if $\xi$ is a L\'evy process, $N$ a
Poisson process and $U=(U_k, k\geq0)$ a sequence of i.i.d. random
variables, $\xi$, $N$ and $U$ being independent, then $Y$ defined by
(\ref{eq444}) is a multiplicative L\'evy process. It is easy to see
that a multiplicative L\'evy process is a symmetric Lamperti--Kiu process.

The reason in Chybiryakov \cite{Chybiryakov} to study the class of multiplicative
L\'evy processes was to establish a Lamperti type representation for
real valued processes that fulfill the scaling property given in the
following definition. A strong Markov family $\{X^{(x)} = (X, \mathbb
{P}_x), {x\in\mathbb{R}^\ast}\}$ with state space $\mathbb{R}^\ast$,
is self-similar of index $\alpha>0$ in the sense of Chybiryakov \cite
{Chybiryakov}, if for all $c\neq0$,
%
\begin{equation}
\label{eq555} \bigl\{(cX_{|c|^{-\alpha} t}, t\geq0), \mathbb{P}_{x}
\bigr\} \overset {\mathcal{L}}= \bigl\{(X_t, t\geq0),
\mathbb{P}_{
cx} \bigr\}
\end{equation}
for all $x\in\mathbb{R}^\ast$. The Lamperti type representation given
in Chybiryakov \cite{Chybiryakov} establishes that for such a self-similar process
$X^{(x)}$, the process $Y$, defined by
\[
Y_t = x^{-1}X_{\nu^{(x)}(t)}^{(x)},\qquad  t\geq0,
\]
where
\[
\nu^{(x)}(t) = \inf \biggl\{s>0\dvt \int_0^s
\bigl|X_u^{(x)}\bigr|^{-\alpha} \,\mathrm{d}u > t \biggr\},\qquad  t\geq0,
\]
with $\inf\{\varnothing\}=\infty$, is a multiplicative L\'evy process.
Conversely, let $Y$ be a multiplicative L\'evy process, and
\[
\mathcal{E}_t = \xi_t + \sum
_{k=1}^{N_t} U_k + \mathrm{i}\ppi N_t,\qquad  t
\geq0,
\]
where $\xi$, $N$ and $(U_k, k\geq0)$ are as in (\ref{eq444}), so that
$Y_t = \exp\{\mathcal{E}_t\}$, $t\geq0$. For $x\in\mathbb{R}^\ast$,
define $X^{(x)}$ by
\[
X_t^{(x)} = xY_{\tau(t x^{-\alpha})},\qquad  t\geq0,
\]
where
\[
\tau(t) = \inf \biggl\{u>0\dvt \int_0^u \bigl|\exp
\{ \alpha\mathcal{E}_u \}\bigr| \,\mathrm{d}u > t \biggr\},\qquad t\geq0,
\]
with $\inf\{\varnothing\}=\infty$. Then $X^{(x)}$ is a $\mathbb{R}^\ast
$-valued self-similar Markov process in the sense of Chybiryakov \cite
{Chybiryakov}, which is recalled in (\ref{eq555}).

It is important to observe that if we take $c=-1$ in (\ref{eq555}), it
is seen that the process $X^{(x)}$ is necessarily a symmetric process
and as a consequence $Y$ is also symmetric. In this work we establish
the analogous description for non-symmetric real valued self-similar
Markov processes.

The remainder of the paper is organized as follows. Section \ref
{rvsMpandLKp} is devoted to some preliminary results about real-valued
self-similar Markov processes. In Section \ref{constofLKp}, we
construct the underlying process in Lamperti's representation and
establish the result that all Lamperti--Kiu processes can be written
this way. Lamperti's representation is given and the infinitesimal
generator of Lamperti--Kiu processes is computed in this section.
Section \ref{proofs} is devoted to prove the main results. In Section
\ref{examples}, we provide two examples where it is possible to compute
explicitly the characteristics of the Lamperti--Kiu process: the $\alpha
$-stable process and the $\alpha$-stable process conditioned to avoid
zero.\looseness=1

\section{Preliminaries and main results}
\subsection{Real-valued self-similar Markov processes and description
of Lamperti--Kiu processes} \label{rvsMpandLKp}

In this section, we will prove some additional properties of
real-valued self-similar Markov processes, in order to characterize
them as time changed Lamperti--Kiu processes.

Let $X$ be a real-valued self-similar Markov process. Let $H_n$ be the
$n$th change of sign of the process $X$, that is,
\[
H_0 = 0, \qquad H_{n} = \inf \{ t>H_{n-1}\dvt
X_tX_{t-} < 0 \},\qquad  n\geq1.
\]
Note that
%
\begin{eqnarray}\label{eq666}
\nonumber
H_1 ( X ) &=& \inf \{ t>0\dvt X_t
X_{t-} < 0 \}
\\
&=& |x|^{\alpha} \inf \bigl\{ |x|^{-\alpha} t>0\dvt \bigl(
|x|^{-1} X_{|x|^{\alpha}|x|^{-\alpha}t} \bigr) \bigl( |x|^{-1}
X_{|x|^{\alpha
}(|x|^{-\alpha}t)-} \bigr) < 0 \bigr\}
\\
&= &|x|^{\alpha} H_1 \bigl( |x|^{-1}
X_{|x|^{\alpha}\cdot} \bigr).\nonumber
\end{eqnarray}
Hence, by the self-similarity property, for $x\in\mathbb{R}^\ast$, it
holds that $\mathbb{P}_x(H_1<\infty)=\mathbb{P}_{\sign(x)}(H_1<\infty
)$. Furthermore, proceeding as in the proof of Lemma 2.5 in Lamperti \cite
{Lamperti72}, it is verified that for each $x\in\mathbb{R}^\ast$,
either $\mathbb{P}_x(H_1<\infty) = 1$ or $\mathbb{P}_x(H_1<\infty) =
0$. The latter and former facts allow us to conclude that there are
four mutually exclusive cases, namely,
\begin{enumerate}
\item[{C.1}] $\mathbb{P}_x(H_1<\infty) = 1$, $\forall  x>0$
and $\mathbb{P}_x(H_1=\infty) = 1$, $\forall  x<0$;

\item[{C.2}] $\mathbb{P}_x(H_1<\infty) = 1$, $\forall  x<0$
and $\mathbb{P}_x(H_1=\infty) = 1$, $\forall  x>0$;

\item[{C.3}] $\mathbb{P}_x(H_1=\infty) = 1$, $\forall  x\in
\mathbb{R}^\ast$;

\item[{C.4}] $\mathbb{P}_x(H_1<\infty) = 1$, $\forall  x\in
\mathbb{R}^\ast$.
\end{enumerate}
In the case {C.1}, if the process $X$ starts at a negative
point, then $\{(-X_t\mathbf{1}_{\{t<T\}}, t\geq0), \mathbb{P}_x\}_{x<0}$
behaves as a positive self-similar Markov process, which have
already been characterized by Lamperti. Now, if the process starts at a
positive point, it can be deduced from Lamperti's representation
(further details are given in the forthcoming Theorem \ref{theo4}(i))
that the process $X$ behaves as a time changed L\'evy process
until it changes of sign, and when this occurs, by the strong Markov
property, its behaviour is that of $X$ issued from a negative point.
The case {C.2} is similar to the first one. For the case
{C.3}, depending on the starting point, $X$ or $-X$ is a positive
self-similar Markov process, again we fall in a known case. In summary,
the Lamperti representation for the cases {C.1}--{C.3} can
be obtained from the Theorem \ref{theo4}{(i)} and the Lamperti
representation for the positive self-similar Markov processes. Thus, we
are only interested in the case {C.4}, where the process $X$
a.s. has at least two changes of sign (and by the strong Markov
property infinitely many changes of sign). For this case, we have the
following proposition.

\begin{proposition} \label{prop2}
If $\mathbb{P}_x(H_1<\infty) = 1$, for all $x\in\mathbb{R}^\ast$, then
the sequence of stopping times $(H_n, n\geq0)$ converges to the first
hitting time of zero $T$, $\mathbb{P}_x$-a.s., for all $x\in\mathbb
{R}^\ast$.
\end{proposition}

The proof of this result will be given in Section \ref{proofs}. We can
see that under the condition of Proposition \ref{prop2}, if $X$ is
killed at $T$, then $X$ has an infinite number of changes of sign
before it dies. Moreover, if $T$ is finite, then $X$ reaches zero at
time $T$ continuously from the left.

The result in Proposition \ref{prop2} is well known in the case where
$X$ is an $\alpha$-stable process and $X$ is not a subordinator. In
that case, if $\alpha\in(0,1]$, $T=\infty$ a.s., while if $\alpha\in
(1,2]$, with probability one, $T<\infty$ and $X$ makes infinitely many
jumps before reaching zero. This process and its Lamperti
representation will be studied in Section \ref{alpha-stablekilledatzero}.

Hereafter, we assume that $\mathbb{P}_x(H_1<\infty) = 1$, for all $x\in
\mathbb{R}^\ast$. Then, for every $n\geq0$, the process $(\mathcal
{X}_t^{(n)}, t\geq0)$ given by
%
\begin{equation}
\label{eq777} \mathcal{X}_t^{(n)} = \frac{ X_{H_n+|X_{H_n}|^\alpha t} }{ |X_{H_n}| },\qquad  0
\leq t< |X_{H_n}|^{-\alpha}(H_{n+1}-H_n),
\end{equation}
is well defined. We call the random variable $X_{H_n}$ an overshoot or
undershoot when $X_{H_n-}<0$ and $X_{H_n}>0$ or $X_{H_n-}>0$ and
$X_{H_n}<0$, respectively. The random variable $X_{H_n-}$ is called the
jump height before crossing of the $x$-axis. The case $X_{H_n-}<0$
means that the change of sign at time $H_n$ is from a negative to a
positive value. Now, we define the sequence of random variables $(J_n,
n\geq0)$ given by the quotient
%
\begin{equation}
\label{eq888} J_n = \frac{X_{H_{n+1}}}{X_{H_{n+1}-}},\qquad n\geq0.
\end{equation}
These random objects satisfy the following properties.

\begin{theorem} \label{theo4}
Let $\{ X^{(x)}=(X, \mathbb{P}_x), {x\in\mathbb{R}^\ast}\}$ be a
family of real-valued self-similar Markov processes of index $\alpha
>0$, such that $\mathbb{P}_x(H_1<\infty) = 1$, for all $x\in\mathbb
{R}^\ast$. Then:
\begin{longlist}[(iii)]
\item[(i)] The paths between sign changes, $(\mathcal{X}^{(n)}, n\geq
0)$, as defined in (\ref{eq777}), are independent under $\mathbb{P}_x$,
for $x\in\mathbb{R}^\ast$. Furthermore, for all $n\geq0$,
%
\begin{equation}
\label{eq999} \bigl\{ \bigl( \mathcal{X}_{t}^{(n)}, 0\leq
t<|X_{H_n}|^{-\alpha
}(H_{n+1}-H_n) \bigr),
\mathbb{P}_x \bigr\} \overset{\mathcal{L}} = \bigl\{ (
X_{t}, 0\leq t<H_{1} ), \mathbb{P}_{\sign
(x)(-1)^n} \bigr
\}.
\end{equation}
Hence, they are time changed L\'evy processes killed at an exponential time.

\item[(ii)] The random variables $J_n, n\geq0$, as defined in (\ref
{eq888}), are independent under $\mathbb{P}_x$, for $x\in\mathbb{R}^\ast
$ and for $n\geq0$, the identity
%
\begin{equation}
\label{eqaux4} \{ J_n, \mathbb{P}_x \} \overset{
\mathcal{L}}= \{ J_0 , \mathbb{P}_{\sign(x)(-1)^n} \},
\end{equation}
holds.

\item[(iii)] For every $n\geq0$, the process $\mathcal{X}^{(n)}$ and
the random variable $J_n$ are independent, under~$\mathbb{P}_x$, for
$x\in\mathbb{R}^\ast$.
\end{longlist}
\end{theorem}

From (\ref{eq999}), we can see that only two independent L\'evy
processes killed at an exponential time are involved in the Lamperti
representation. In the same way, from (\ref{eqaux4}), only two
independent real random variables represent the quotient between
overshoots (undershoots) and jump height before crossing of the
$x$-axis.\vadjust{\goodbreak}  Furthermore, by {(iii)} all these random objects are
independent. The latter theorem is at the heart of our motivation to
construct the Lamperti--Kiu processes in the next section.

\subsection{Construction of Lamperti--Kiu processes} \label{constofLKp}

In this section, we give a generalization of time changed exponentials
of L\'evy processes as well as of the processes which are defined in
(\ref{eq444}). We will see that all Lamperti--Kiu processes can be
constructed as this generalization of (\ref{eq444}).

Let $\xi^+$, $\xi^-$ be real valued L\'evy processes; $\zeta^+, \zeta^-$ exponential random variables with parameters $q^+, q^-$,
respectively, and $ U^+, U^-$ real valued random variables. Let $(\xi^{+,k}, k\geq0)$, $(\xi^{-,k}, k\geq0)$, $(\zeta^{+,k}, k\geq0)$,
$(\zeta^{-,k}, k\geq0)$, $(U^{+,k}, k\geq0)$, $(U^{-,k}, k\geq0)$ be
independent sequences of i.i.d. random variables such that
\begin{eqnarray*}
\xi^{+,0}& \overset{\mathrm{Law}} {=} &\xi^+,\qquad \xi^{-,0} \overset
{\mathrm{Law}} {=} \xi^-, \qquad\zeta^{+,0} \overset{\mathrm{Law}} {=}
\zeta^+,\qquad \zeta^{-,0} \overset{\mathrm{Law}} {=} \zeta^-,\\
 U^{+,0}
&\overset{\mathrm{Law}} {=}& U^+, \qquad U^{-,0} \overset{\mathrm {Law}} {=} U^-.
\end{eqnarray*}
For every $x\in\mathbb{R}^\ast$ fixed, we consider the sequence $((\xi^{(x,k)}, \zeta^{(x,k)}, U^{(x,k)}), k\geq0)$, where for \mbox{$k\geq0$},
\[
\bigl(\xi^{(x,k)}, \zeta^{(x,k)}, U^{(x,k)} \bigr) = \cases{
\bigl(\xi^{+,k},
\zeta^{+,k}, U^{+,k} \bigr), &\quad  $\mbox{if $\sign(x) (-1)^k
= 1$,}$
\vspace*{2pt}\cr
 \bigl(\xi^{-,k}, \zeta^{-,k}, U^{-,k} \bigr), & \quad $\mbox{if $\sign(x) (-1)^k =
 -1.$}$}
\]
Let $(T_n^{(x)}, n\geq0)$ be the sequence defined by
\[
T_0^{(x)} = 0, \qquad T_n^{(x)} = \sum
_{k=0}^{n-1} \zeta^{(x,k)},\qquad  n\geq1,
\]
and $(N_t^{(x)}, t\geq0)$ be the alternating renewal type process:
\[
N_t^{(x)} = \max \bigl\{n\geq0\dvt T_n^{(x)}
\leq t \bigr\},\qquad  t\geq0.
\]
For notational convenience, we write
\[
\sigma_t^{(x)} = t-T_{N_t^{(x)}}^{(x)},\qquad
\xi_{ \sigma_t }^{(x)} = \xi_{ \sigma_t^{(x)} }^{(x,N_t^{(x)})},\qquad
\xi_{\zeta}^{(x,k)} = \xi_{\zeta^{(x,k)}}^{(x,k)}.
\]
Finally, we define the process $Y^{(x)} = (Y_t^{(x)}, t\geq0)$ by
%
\begin{equation}
\label{eq1010} Y_t^{(x)} = x\exp \bigl\{
\mathcal{E}_t^{(x)} \bigr\},\qquad  t\geq0,
\end{equation}
where
\[
\mathcal{E}_t^{(x)} = \xi_{ \sigma_t }^{(x)} +
\sum_{k=0}^{N_t^{(x)}-1} \bigl( \xi_{\zeta}^{(x,k)}
+ U^{(x,k)} \bigr) + \mathrm{i}\ppi N_t^{(x)},\qquad  t\geq0.
\]

\begin{remark}
Observe that the process $Y^{(x)}$ is a generalization of
multiplicative L\'evy processes. For, take $(\xi^+, U^+, \zeta^+)
\overset{\mathcal{L}}= (\xi^-, U^-, \zeta^-)$ it is seen that $Y^{(x)}$
is a multiplicative L\'evy process, as it has been defined in Chybiryakov \cite
{Chybiryakov}. Moreover, if $q^+ = 0$ and $q^->0$, then for $x>0$,
$Y^{(x)}$ does not jump to the negative axis and $Y^{(x)}$ is the
exponential of a L\'evy process, which appears in the Lamperti
representation for positive self-similar Markov processes.
\end{remark}

The following theorem is the main result of this paper. The first part
states that $Y^{(x)}$ is a Lamperti--Kiu process, the second and third
parts are the generalization of the Lamperti representation.

\begin{theorem} \label{theo7}
Let $Y^{(x)}$ be the process defined in (\ref{eq1010}). Then
\begin{longlist}[(iii)]
\item[(i)] the process $Y^{(x)}$ is Fellerian in $\mathbb{R}^\ast$ and
satisfies (\ref{eq333}). Furthermore, for any finite stopping time
$\mathbf{T}$:
\[
\bigl( \bigl(Y_{\mathbf{T}}^{(x)} \bigr)^{-1}
Y_{\mathbf{T}+s}^{(x)}, s\geq0 \bigr) \overset {\mathcal{L}}= \bigl(
\exp \bigl\{ \widetilde{\mathcal{E}}_s^{(\sign(Y_\mathbf
{T}^{(x)}))} \bigr\}, s\geq0
\bigr),
\]
where $\widetilde{\mathcal{E}}^{(\cdot)}$ is a copy of $\mathcal
{E}^{(\cdot)}$ which is independent of $(\mathcal{E}_u^{(\cdot)}, 0\leq
u\leq\mathbf{T})$.

\item[(ii)] Let $\{ X^{(x)}=(X, \mathbb{P}_x), {x\in\mathbb{R}^\ast}\}
$ be a family of real-valued self-similar Markov processes of index
$\alpha>0$ such that $\mathbb{P}_x(H_1<\infty)=1$, for all $x\in\mathbb
{R}^\ast$. For every $x\in\mathbb{R}^\ast$ define the process $\mathcal
{Y}^{(x)}$ by
\[
\mathcal{Y}_t^{(x)} = X_{\nu^{(x)}(t)}^{(x)},\qquad  t
\geq0,
\]
where
\[
\nu^{(x)}(t) = \inf \biggl\{s>0\dvt \int_0^s
\bigl|X_u^{(x)}\bigr|^{-\alpha} \,\mathrm{d}u > t \biggr\}.
\]
Then $\mathcal{Y}^{(x)}$ may be decomposed as in (\ref{eq1010}).
Moreover, every Lamperti--Kiu process can be constructed as explained in
(\ref{eq1010}).

\item[(iii)] Conversely, let $(Y^{(x)})_{x\in\mathbb{R}^\ast}$ be a
family of processes as constructed in (\ref{eq1010}) and consider the
processes $ X^{(x)}$ given by
\[
X_t^{(x)} = Y_{\tau(t|x|^{-\alpha})}^{(x)},\qquad  t\geq0,
\]
where
\[
\tau(t) = \inf \biggl\{s>0\dvt \int_0^s \bigl|\exp
\bigl\{\alpha\mathcal{E}_u^{(x)} \bigr\}\bigr| \,\mathrm{d}u > t \biggr\},\qquad
t<T
\]
for some $\alpha>0$. Then $(X^{(x)})_{x\in\mathbb{R}^\ast} $ is a
family of real-valued self-similar Markov processes of index $\alpha>0$.
\end{longlist}
\end{theorem}

From now on, we denote a Lamperti--Kiu process by $Y$. Now, we obtain an
expression for the infinitesimal generator of $Y$, that will be used in
the examples.

\begin{proposition} \label{prop8}
Let $\mathcal{K}$ be the infinitesimal generator of $Y$. Let $\mathcal
{A}^+$, $\mathcal{A}^-$ be the infinitesimal generators of $\xi^+$, $\xi^-$,
respectively. Let $f$ be a bounded continuous function such that
$f(0) = 0$ and $(f\circ\exp) \in\mathcal{D}_{\mathcal{A}^{+}}$ and
$(f\circ-\exp) \in\mathcal{D}_{\mathcal{A}^{-}}$. Then, for every
$x\in\mathbb{R}^\ast$,
%
\begin{equation}
\label{eq1111} \mathcal{K}f(x) = \mathcal{A}^{\sign(x)} \bigl(f\circ\sign(x)
\exp \bigr) \bigl(\log|x|\bigr) + q^{\sign(x)} \bigl( \mathbf{E} \bigl[f \bigl(-x\exp
\bigl\{U^{\sign(x)} \bigr\} \bigr) - f(x) \bigr] \bigr).
\end{equation}
\end{proposition}

With the help of the latter proposition, we can give the infinitesimal
generator of Y in terms of the parameters of the L\'evy processes
$\xi^+$ and $\xi^-$ as follows. Recall that the characteristic exponent of
the L\'evy process $\xi^{\pm}$ can be written as
\[
\psi^{\pm}(\lambda) = a^{\pm}\mathrm{i}\lambda- \frac{[\sigma^{\pm
}]^2}{2}
\lambda^2 + \int_{\mathbb{R}} \bigl[\mathrm{e}^{\mathrm{i}\lambda y} - 1
- \mathrm{i}\lambda l(y) \bigr] \pi^{\pm}(\mathrm{d}y),\qquad \lambda\in\mathbb{R},
\]
where $a^{\pm}\in\mathbb{R}$, $\sigma>0$, $l(\cdot)$ is a fixed
continuous bounded function such that $l(y)\sim y$ as $y\rightarrow0$
and $\pi^{\pm}$ is the L\'evy measure of the process $\xi^{\pm}$, which
satisfies $\pi^{\pm}(\{0\})=0$ and $\int_{\mathbb{R}} (1\wedge x^2)\pi^{\pm}(\mathrm{d}x)<\infty$.
Furthermore, the choice of the function $l$ is
arbitrary and the coefficient $a^{\pm}$ is the only one which depends
on this choice (see Remark 8.4 in Sato \cite{Sato}). Later in the examples,
we will choose conveniently this function. Hence, the infinitesimal
generator of the L\'evy process $\xi^{\pm}$ can be expressed as
\[
\mathcal{A}^{\pm}f(x) = a^{\pm}f'(x) +
\frac{[\sigma^{\pm}]^2}{2}f''(x) + \int_{\mathbb{R}}
\bigl[f(x+y) - f(x) - f'(x)l(y) \bigr] \pi^{\pm}(\mathrm{d}y),\qquad  f\in
\mathcal{D}_{\mathcal{A}^{\pm}}.
\]
Then using the expression of $\mathcal{A}^{\pm}$ and (\ref{eq1111}), we
find for $x\in\mathbb{R}^\ast$,
%
\begin{eqnarray}
\label{eq1212} \mathcal{K}f(x) &=& b^{\sign(x)}xf'(x) +
\frac{[\sigma^{\sign
(x)}]^2}{2}x^2f''(x)\nonumber
\\
&&{} + q^{\sign(x)} \bigl\{ \mathbf{E} \bigl[f \bigl(-x\exp \bigl
\{U^{\sign(x)} \bigr\} \bigr) - f(x) \bigr] \bigr\}
\\
&&{} + \int_{\mathbb{R}^+} \bigl[f(xu) - f(x) - xf'(x)l(
\log u) \bigr] \Theta^{\sign(x)}(\mathrm{d}u),\nonumber
\end{eqnarray}
where $b^{\sign(x)} = a^{\sign(x)}+[\sigma^{\sign(x)}]^2/2$,
$\Theta^{\sign(x)}(\mathrm{d}u) = \pi^{\sign(x)}(\mathrm{d}u)\circ\log u$. Hence, by
Volkonskii's theorem, the generator $\widetilde{\mathcal{K}}$ of the
time changed process $Y_{\tau}$ is given by $\widetilde{\mathcal
{K}}f(x) = |x|^{-\alpha}\mathcal{K}f(x)$, for $x\in\mathbb{R}^\ast$.
Hence, knowing that the infinitesimal generator of $Y$ is given by (\ref
{eq1212}) it is possible to identify the infinitesimal generator of the
self-similar Markov process $X$ and conversely.

\section{Proofs} \label{proofs}\vspace*{-9pt}

\begin{pf*}{Proof of Proposition \ref{prop2}}
The strong Markov property implies $\mathbb{P}_x(H_n<\infty, \forall
n\geq0) = 1$. Thus, $(H_n, n\geq0)$ is a strictly increasing sequence
of stopping times satisfying $H_n\leq T$, for all $n\geq0$. Let $H$ be
the limit of this sequence, then $H\leq T$. If $H=\infty$, then clearly
$T=\infty$ and $H=T$. On the other hand, if $H<\infty$, then on the set
$\{H<T\}$, it is possible to define the process $X^H = (X_{H+t}\mathbf
{1}_{\{t<T-H\}}, t\geq0)$. This process has no change of sign, and by
the strong Markov property, for all $y\in\mathbb{R}^\ast$,
conditionally on $X_H=y$, $X^H$ has the same distribution as $X$ under
$\mathbb{P}_y$. This contradicts the fact that $X$ has at least one
change of sign. Therefore, $H=T$, a.s.
\end{pf*}

\begin{pf*}{Proof of Theorem \ref{theo4}}
For $t\geq0$, we denote by $\theta_t\dvt \Omega\rightarrow\Omega$ the
shift operator, that is, for $\omega\in\Omega$, $\theta_t\omega(s) =
\omega(t+s)$, $s\geq0$.
\begin{longlist}[(iii)]
\item[(i)] Let $F$ be a bounded and measurable functional. From (\ref
{eq666}) and the self-similarity property, it follows
\[
\mathbb{E}_x \biggl[ F \biggl( \frac{ X_{|X_{0}|^{\alpha} t} }{ |X_0| }, 0\leq
t<|X_{0}|^{-\alpha}H_1 \biggr) \biggr] =
\mathbb{E}_{\sign(x)} \bigl[ F ( X_{t}, 0\leq t<H_1 )
\bigr]
\]
for $x\in\mathbb{R}^\ast$. Moreover, $\sign(X_{H_n}) = \sign
(x)(-1)^n$, $\mathbb{P}_x$-a.s. These two facts and the strong Markov
property are sufficient to complete the proof. Indeed, for $\mathcal
{X}^{(0)}, \ldots, \mathcal{X}^{(n)}$ as defined in (\ref{eq777}) and
for all $F_0, \ldots, F_n$ bounded and measurable functionals, we have
\begin{eqnarray*}
\mathbb{E}_x \Biggl[ \prod_{k=0}^n
F_k \bigl( \mathcal{X}^{(k)} \bigr) \Biggr] &=&
\mathbb{E}_x \Biggl[ \prod_{k=0}^{n-1}
F_k \bigl( \mathcal{X}^{(k)} \bigr) \mathbb{E}_{X_{H_n}}
\biggl[ F_n \biggl( \frac{ X_{|X_{0}|^{\alpha
} t} }{ |X_0| }, 0\leq t<|X_{0}|^{-\alpha}H_1
\biggr) \biggr] \Biggr]
\\
&= &\mathbb{E}_x \Biggl[ \prod_{k=0}^{n-1}
F_k \bigl( \mathcal{X}^{(k)} \bigr) \mathbb{E}_{\sign(x)(-1)^n}
\bigl[ F_n ( X_t, 0\leq t<H_1 ) \bigr]
\Biggr]
\\
&=& \mathbb{E}_x \Biggl[ \prod_{k=0}^{n-1}
F_k \bigl( \mathcal{X}^{(k)} \bigr) \Biggr]
\mathbb{E}_{\sign(x)(-1)^n} \bigl[ F_n ( X_t, 0\leq
t<H_1 ) \bigr],
\end{eqnarray*}
where the strong Markov and self-similarity properties were used to
obtain the first and second equality, respectively. Now, taking $F_0 =
\cdots= F_{n-1} \equiv1$, we have
\[
\mathbb{E}_x \bigl[ F_n \bigl( \mathcal{X}_{t}^{(n)},
0\leq t<|X_{H_n}|^{-\alpha}(H_{n+1}-H_n)
\bigr) \bigr] = \mathbb{E}_{\sign
(x)(-1)^n} \bigl[ F_n (
X_t, 0\leq t<H_1 ) \bigr].
\]
This proves (\ref{eq999}). In addition
\[
\mathbb{E}_x \Biggl[ \prod_{k=0}^n
F_k \bigl( \mathcal{X}^{(k)} \bigr) \Biggr] =
\mathbb{E}_x \Biggl[ \prod_{k=0}^{n-1}
F_k \bigl( \mathcal {X}^{(k)} \bigr) \Biggr]
\mathbb{E}_x \bigl[ F_n \bigl( \mathcal
{X}^{(n)} \bigr) \bigr] .
\]
This proves the independence in the sequence $\{(\mathcal{X}_t^{(n)},
0\leq t<|X_{H_n}|^{-\alpha}(H_{n+1}-H_n)), n\geq0\}$ under $\mathbb{P}_x$.

\item[(ii)] From (\ref{eq666}) and the self-similarity property, we
derive that
%
\begin{equation}
\label{eq1313} \mathbb{E}_x \biggl[ f \biggl( \frac{ X_{H_1} }{ X_{H_1-} }
\biggr) \biggr] = \mathbb{E}_{\sign(x)} \biggl[ f \biggl( \frac{ X_{H_1} }{
X_{H_1-} }
\biggr) \biggr]
\end{equation}
for all $x\in\mathbb{R}^\ast$, and $f$ bounded Borel function. Now,
let $f_0, \ldots, f_n$ bounded Borel functions. Proceeding as in \emph
{(i)}, using (\ref{eq1313}) and the strong Markov property, we obtain
\begin{eqnarray*}
&&\mathbb{E}_x \Biggl[ \prod_{k=0}^n
f_k \biggl( \frac{ X_{H_{k+1}}
}{X_{H_{k+1}-}} \biggr) \Biggr] \\
&&\quad=
\mathbb{E}_x \Biggl[ \prod_{k=0}^{n-1}
f_k \biggl( \frac{ X_{H_{k+1}} }{X_{H_{k+1}-}} \biggr) \Biggr] \mathbb
{E}_{\sign(x)(-1)^n} \biggl[ f_n \biggl( \frac{ X_{H_{1}} }{X_{H_{1}-}} \biggr)
\biggr] .
\end{eqnarray*}
The conclusion follows as in {(i)}.

\item[(iii)] By the strong Markov property, {(i)} and {(ii)},
it is sufficient to prove the case $n=0$. For $k\geq1$, let $f\dvt {\mathbb
{R}^{\ast}}^k \rightarrow\mathbb{R}$, $g\dvt \mathbb{R}^{\ast} \rightarrow
\mathbb{R}$ be two Borel functions, and $0<s_1<\cdots<s_k$. We note
the following identity $\frac{ X_{H_1} }{ X_{H_1-} } \circ\theta_{s_k}
= \frac{ X_{H_1} }{ X_{H_1-} }$, on $\{ s_k< H_1 \}$. Hence, by the
Markov property and (\ref{eq1313}), we have
\begin{eqnarray*}
&&\mathbb{E}_x \biggl[ f(X_{s_1}, \ldots,
X_{s_k}) g \biggl( \frac{ X_{H_1}
}{ X_{H_1-} } \biggr);s_k<H_1
\biggr]
\\
&&\quad= \mathbb{E}_x \biggl[ f(X_{s_1}, \ldots,
X_{s_k}) \mathbb{E}_{X_{s_k}} \biggl[ g \biggl(
\frac{ X_{H_1} }{ X_{H_1-} } \biggr) \biggr]; s_k<H_1 \biggr]
\\
&&\quad= \mathbb{E}_x \bigl[ f(X_{s_1}, \ldots,
X_{s_k}); s_k<H_1 \bigr] \mathbb{E}_{\sign(x)}
\biggl[ g \biggl( \frac{ X_{H_1} }{ X_{H_1-} } \biggr) \biggr]
\\
&&\quad= \mathbb{E}_x \bigl[ f(X_{s_1}, \ldots,
X_{s_k}); s_k<H_1 \bigr] \mathbb{E}_{x}
\biggl[ g \biggl( \frac{ X_{H_1} }{ X_{H_1-} } \biggr) \biggr].
\end{eqnarray*}
This ends the proof.
\end{longlist}
\end{pf*}

In order to prove Theorem \ref{theo7}, we first prove the following
lemma. This lemma is a consequence of the lack-of-memory property of
the exponential distribution and the properties of the random objects
which define $Y^{(x)}$. Before we state it, we define the following
process. For $x\in\mathbb{R}^\ast$, let $Z^{(x)}$ be the sign process
of $Y^{(x)}$, that is, $Z_t^{(x)} = \sign(Y_t^{(x)})$, $t\geq0$. Note
that $Z^{(x)}$ is a continuous time Markov chain with state space $\{
-1,1\}$, starting point $\sign(x)$ and transition semigroup $\mathrm{e}^{t\mathbf
{Q}}$, where
\[
\mathbf{Q} = \pmatrix{ -q^- & q^-
\vspace*{2pt}\cr
q^+ & -q^+ }.
\]
Furthermore, since the law of $Z^{(x)}$ is determined by $\mathbf{Q}$
(hence by $\zeta^+, \zeta^-$), then the process $Z^{(x)}$ is
independent of $((\xi^{(x,k)}, U^{(x,k)}), k\geq0)$.

\begin{lemma} \label{lemma10}
Let $n, m$ be positive integers and $s,t$ be positive real numbers. We
have the following properties
\begin{longlist}[(a)]
\item[(a)] Conditionally on $T_n^{(x)}\leq t < T_{n+1}^{(x)}$, the
random variable $T_{n+1}^{(x)}-t$ has an exponential distribution with
parameter $q^{(x,n)}$, where $q^{(x,n)}$ equals $q^+$ if $\sign
(x)(-1)^n=1$ and $q^-$ otherwise. Furthermore,
\[
\xi_{\zeta}^{(x,n)} - \xi_{t-T_n^{(x)}}^{(x,n)}
\overset{\mathcal{L}}= \widetilde{\xi}_{\widetilde{\zeta}}^{(Z_t^{(x)},0)},
\]
where $(\widetilde{\xi}^{(\cdot,0)}, \widetilde{\zeta}^{(\cdot, 0)})$
are independent of $(\xi^{(\cdot, k)}, \zeta^{(\cdot, k)}, 0\leq k< n
)$ and with the same distribution as $(\xi^{(\cdot,0)}, \zeta^{(\cdot, 0)})$.

\item[(b)] Conditionally on $T_n^{(x)}\leq t < T_{n+1}^{(x)}$,
$T_{n+m}^{(x)}\leq t+s < T_{n+m+1}^{(x)}$ the distribution of $\xi_{t+s-T_{n+m}^{(x)}}^{(x,n+m)}$ is the same as the distribution of
$\widetilde{\xi}_{s-\widetilde{T}_m^{(Z_t^{(x)})}}^{(Z_t^{(x)},m)}$
conditionally on $\widetilde{T}_{m}^{(Z_t^{(x)})} \leq s < \widetilde
{T}_{m+1}^{(Z_t^{(x)})}$, that is,
\begin{eqnarray*}
&&\mathbf{P} \bigl( \xi_{t+s - T_{n+m}^{(x)}}^{(x,n+m)} \in \mathrm{d}z \mid
T_{n+m}^{(x)}\leq t+s < T_{n+m+1}^{(x)},
T_n^{(x)}\leq t < T_{n+1}^{(x)} \bigr)
\\
&&\quad= \mathbf{P} \bigl( \widetilde{\xi}_{s - \widetilde
{T}_{m}^{(Z_t^{(x)})}}^{(Z_t^{(x)},m)} \in \mathrm{d}z \mid
\widetilde {T}_{m}^{(Z_t^{(x)})}\leq s < \widetilde{T}_{m+1}^{(Z_t^{(x)})}
\bigr),
\end{eqnarray*}
where $(\widetilde{\xi}^{(\cdot,m)}, \widetilde{T}_m^{(\cdot)})$ are
independent of $(\xi^{(\cdot, k)}, T_k^{(\cdot)}, 0\leq k\leq n )$ with
the same distribution as $(\xi^{(\cdot,m)}, T_m^{(\cdot)})$.
\end{longlist}
\end{lemma}

\begin{pf*}{Proof of Lemma \ref{lemma10}}
The first part of (a) follows from the lack-of-memory property of the
exponential distribution. Now, by construction, $(\xi^{(x,n)}, \zeta^{(x,n)}, n\geq0)$ is a sequence of independent random objects which
depends on $x$ only through its sign and $T_{n+m}^{(x)} = T_n^{(x)} +
\sum_{k=0}^{m-1} \zeta^{(x,n+k)}$. Hence, it is always possible to take
$(\widetilde{\xi}^{(\cdot,0)}, \widetilde{\zeta}^{(\cdot, 0)})$ and
$(\widetilde{\xi}^{(\cdot,m)},\widetilde{T}_m^{(\cdot)})$ with the
properties described in (a) and (b), respectively. Thus, it only
remains to prove the equality in distribution in (a) and~(b).

Denote by $f_{T_n^{(x)}}$ the density of the random variable
$T_n^{(x)}$. Simple computations lead to
\begin{eqnarray*}
&&\mathbf{P} \bigl(\xi_{\zeta}^{(x,n)} - \xi_{t-T_n^{(x)}}^{(x,n)}
\in \mathrm{d}z, T_n^{(x)}\leq t < T_{n+1}^{(x)}
\bigr)
\\
&&\quad= \int_0^t \int_{t-u}^\infty
\mathbf{P} \bigl(\xi_{r-(t-u)}^{(x,n)}\in \mathrm{d}z \bigr)
q^{(x,n)}\mathrm{e}^{-q^{(x,n)}r}\,\mathrm{d}r f_{T_{n}^{(x)}}(u)\,\mathrm{d}u
\\
&&\quad= \mathbf{P} \bigl(\xi_{\zeta}^{(x,n)} \in \mathrm{d}z \bigr) \mathbf{P}
\bigl(T_n^{(x)}\leq t < T_{n+1}^{(x)}
\bigr),
\end{eqnarray*}
where the independence and stationarity of the increments of the L\'evy
process $\xi^{(x,n)}$ have been used in the first equality and we made
the change of variables $v=r-(t-u)$ to obtain the second. Hence, the
equality in law of (a) is obtained.

By (a), we have that for all $m\geq0$, conditionally on $T_n^{(x)}\leq
t < T_{n+1}^{(x)}$, the random variable $T_{n+m}^{(x)}-t$ has the same
distribution as $T_m^{(Z_t^{(x)})}$ and it is independent of
$(T_k^{(x)}, 0\leq k\leq n)$. Hence,
\begin{eqnarray*}
&&\mathbf{P} \bigl( \xi_{t+s - T_{n+m}^{(x)}}^{(x,n+m)} \in \mathrm{d}z,
T_{n+m}^{(x)}\leq t+s < T_{n+m+1}^{(x)} \mid
T_n^{(x)}\leq t < T_{n+1}^{(x)} \bigr)
\\
&&\quad= \mathbf{P} \bigl( \xi_{s - \widetilde
{T}_{m}^{(Z_t^{(x)})}}^{(Z_t^{(x)},m)} \in \mathrm{d}z, \widetilde
{T}_{m}^{(Z_t^{(x)})}\leq s < \widetilde{T}_{m+1}^{(Z_t^{(x)})}
\bigr)
\end{eqnarray*}
and
\[
\mathbf{P} \bigl(T_{n+m}^{(x)}\leq t+s < T_{n+m+1}^{(x)}
\mid T_n^{(x)}\leq t < T_{n+1}^{(x)}
\bigr) = \mathbf{P} \bigl(\widetilde{T}_{m}^{(Z_t^{(x)})}\leq s <
\widetilde{T}_{m+1}^{(Z_t^{(x)})} \bigr).
\]
Therefore,
\begin{eqnarray*}
&&\mathbf{P} \bigl( \xi_{t+s - T_{n+m}^{(x)}}^{(x,n+m)} \in \mathrm{d}z \mid
T_{n+m}^{(x)}\leq t+s < T_{n+m+1}^{(x)},
T_n^{(x)}\leq t < T_{n+1}^{(x)} \bigr)
\\
&&\quad= \mathbf{P} \bigl( \xi_{s - \widetilde
{T}_{m}^{(Z_t^{(x)})}}^{(Z_t^{(x)},m)} \in \mathrm{d}z \mid\widetilde
{T}_{m}^{(Z_t^{(x)})}\leq s < \widetilde{T}_{m+1}^{(Z_t^{(x)})}
\bigr).
\end{eqnarray*}
This finishes the proof.
\end{pf*}

\begin{pf*}{Proof of Theorem \ref{theo7}}
{(i)} First, we prove that $Y^{(x)}$ satisfies the property (\ref
{eq333}). We note that the process $\mathcal{E}^{(\cdot)}$ depends on
$x$ only through its sign, then clearly for all $a\in\mathbb{R}^\ast$,
$\mathcal{E}^{(|a|x)} \overset{\mathcal{L}}= \mathcal{E}^{(x)}$. Hence,
we have
\begin{eqnarray*}
\bigl(\sign(a)Y_t^{(|a|x)}, t\geq0 \bigr) &=& \bigl(
\sign(a)|a|x\exp \bigl\{ \mathcal {E}_t^{(|a|x)} \bigr\}, t\geq0
\bigr)
\\
&\stackrel{\mathcal{L}} = & \bigl(ax\exp \bigl\{\mathcal{E}_t^{(x)}
\bigr\}, t\geq0 \bigr)
\\
&=& \bigl(aY_t^{(x)}, t\geq0 \bigr).
\end{eqnarray*}
Therefore, the process $Y^{(x)}$ satisfies the property (\ref{eq333}).

Let $s,t\geq0$, then by Lemma \ref{lemma10}, conditionally on
$T_n^{(x)}\leq t < T_{n+1}^{(x)}$, $T_{n+m}^{(x)}\leq t+s <
T_{n+m+1}^{(x)}$, we have
\begin{eqnarray*}
\frac{ Y_{t+s}^{(x)} }{ Y_t^{(x)} } &=& \exp \Biggl\{ \xi_{t+s-T_{n+m}^{(x)}}^{(x,n+m)} + \sum
_{k=1}^{m-1} \bigl( \xi_{\zeta
}^{(x,n+k)}
+ U^{(x,n+k)} \bigr) + \xi_{\zeta}^{(x,n)} -
\xi_{t-T_n^{(x)}}^{(x,n)} + U^{(x,n)} + \mathrm{i}\ppi m \Biggr\}
\\
&\overset{\mathcal{L}}=& \exp \Biggl\{ \widetilde{\xi}_{s-\widetilde
{T}_{m}^{(Z_t^{(x)})}}^{(Z_t^{(x)},m)}
+ \sum_{k=0}^{m-1} \bigl( \widetilde{
\xi}_{ \widetilde{\zeta} }^{(Z_t^{(x)},k)} + \widetilde {U}^{(Z_t^{(x)},k)} \bigr) + \mathrm{i}\ppi
m \Biggr\}.
\end{eqnarray*}
Hence, for $s,t\geq0$,
%
\begin{equation}
\label{eq1414} \frac{ Y_{t+s}^{(x)} }{ Y_t^{(x)} } \overset{\mathcal{L}}= \exp \bigl\{
\widetilde{ \mathcal{E}}_s^{(Z_t^{(x)})} \bigr\},
\end{equation}
where $\widetilde{\mathcal{E}}^{(\cdot)}$ is a copy of $\mathcal
{E}^{(\cdot)}$ which is independent of $( \mathcal{E}_u^{(\cdot)},
0\leq u\leq t)$. Thus, $Y^{(x)}$ has the Markov property. Furthermore,
\[
\bigl( Y_{t+s}^{(x)}, s\geq0 \bigr) \overset{\mathcal{L}} =
\bigl( \widetilde {Y}_s^{(Y_t^{(x)})}, s\geq0 \bigr),
\]
where $\widetilde{Y}^{(\cdot)}$ is a copy of $Y^{(\cdot)}$ which is
independent of $( Y_u^{(\cdot)}, 0\leq u\leq t)$. This also ensures
that all processes $Y^{(x)}$ have the same semigroup.

Now, we prove that $Y^{(x)}$ is a Feller process on $\mathbb{R}^\ast$.
Let $Q_t$ be the semigroup associated to $Y^{(x)}$. We verify that
$Q_t$ is a Feller semigroup, that is,
\begin{longlist}[(ii)]
\item[(i)] $Q_tf \in C_0(\mathbb{R}^\ast)$, for all $f\in C_0(\mathbb
{R}^\ast),$

\item[(ii)] $\lim_{t\downarrow0}Q_tf(x) = f(x)$, for all $x \in\mathbb
{R}^\ast$.
\end{longlist}
Let $x\in\mathbb{R}^\ast$ be fixed. For all $y\in\mathbb{R}^\ast$
such that $\sign(y) = \sign(x)$, by property (\ref{eq333}), we have
\[
Q_tf(y) = \mathbf{E} \bigl[ f \bigl( Y_t^{(y)}
\bigr) \bigr] = \mathbf {E} \biggl[ f \biggl( \frac{y}{x}
Y_t^{(x)} \biggr) \biggr].
\]
The latter expression and the dominated convergence theorem ensure the
continuity of $Q_tf$ in~$x$. By (\ref{eq333}), $(Y_t^{(x)}, t\geq0)
\stackrel{\mathcal{L}} = (|x|Y_t^{(\sign(x))}, t\geq0)$ for all $x \in
\mathbb{R}^\ast$. Hence,
\[
Q_tf(x) = \mathbf{E} \bigl[ f \bigl( Y_t^{(x)}
\bigr) \bigr] = \mathbf {E} \bigl[ f \bigl( |x|Y_t^{(\sign(x))}
\bigr) \bigr], \qquad x\in\mathbb{R}^\ast.
\]
Using again the dominated convergence theorem, we obtain $\lim_{|x|\rightarrow\infty} Q_tf(x) = 0$. For the last part,
\begin{eqnarray*}
\mathbf{E} \bigl[ f \bigl( Y_t^{(x)} \bigr) \bigr] &=&
\mathbf{E} \bigl[ f \bigl( Y_t^{(x)} \bigr) \rrvert
T_1^{(x)} > t \bigr] \mathbf{P} \bigl( T_1^{(x)}
> t \bigr)
\\
&&{} + \mathbf{E} \bigl[ f \bigl( Y_t^{(x)} \bigr) \rrvert
T_1^{(x)} \leq t \bigr] \mathbf{P} \bigl(
T_1^{(x)} \leq t \bigr).
\end{eqnarray*}
For the first term, we have
\[
\mathbf{E} \bigl[ f \bigl( Y_t^{(x)} \bigr) \rrvert
T_1^{(x)} > t \bigr] \mathbf{P} \bigl( T_1^{(x)}
> t \bigr) = \mathbf{E} \bigl[ f \bigl( x\exp \bigl\{ \xi_t^{\sign(x)}
\bigr\} \bigr) \rrvert \zeta^{\sign(x)} > t \bigr] \mathrm{e}^{-q^{\sign(x)}t}.
\]
Letting $t\rightarrow0$, the last expression converges to $f(x)$ by
the right continuity of $\xi^{\sign(x)}$. Thus, it only remains to
prove that the second term converges to zero as $t$ tends to zero.
Since $f$ is bounded,
\[
\bigl\llvert \mathbf{E} \bigl[ f \bigl( Y_t^{(x)} \bigr)
\rrvert T_1^{(x)} \leq t \bigr] \mathbf{P} \bigl(
T_1^{(x)} \leq t \bigr) \bigr\rrvert \leq C \bigl(
1-\mathrm{e}^{-q^{\sign(x)}t} \bigr)
\]
for some positive constant $C$. Again, letting $t\rightarrow0$ we
obtain the desired result.

The strong Markov property of $Y^{(x)}$ follows from the standard fact
that any Feller process is a strong Markov process.

{(ii)} First, note that $\nu^{(x)}(t)$ satisfies
%
\begin{equation}
\label{eq1515} \nu^{(x)}(t) = \int_0^t
\bigl|\mathcal{Y}_s^{(x)}\bigr|^{\alpha}\,\mathrm{d}s, \qquad t\geq0.
\end{equation}
Indeed, if
\[
\tau^{(x)}(t) = \int_0^{t|x|^{\alpha}}
\bigl|X_s^{(x)}\bigr|^{-\alpha}\,\mathrm{d}s,
\]
then, since $\tau^{(x)}( \nu^{(x)}(t)|x|^{-\alpha} ) = t$, it follows
$\,\mathrm{d}\nu^{(x)}(t)/\,\mathrm{d}t = 1/|X_{\nu^{(x)}(t)}^{(x)}|^{-\alpha} = |\mathcal
{Y}_t^{(x)}|^{\alpha}$.

Now, we claim the following: for every $x\in\mathbb{R}^\ast$ and
$n\geq0$, there exists a L\'evy process $\xi^{(x,n)}$ independent of
$(X_s^{(x)}, 0\leq s\leq H_n^{(x)})$ such that,
%
\begin{equation}
\label{eq1616} X_{H_n+t}^{(x)} = X_{H_n}^{(x)}
\exp \bigl\{ \xi_{\tau
^{(x,n)}(t|X_{H_n}^{(x)}|^{-\alpha})}^{(x,n)} \bigr\}, \qquad 0\leq t<
H_{n+1}^{(x)}-H_n^{(x)},
\end{equation}
where
%
\begin{equation}
\label{eq1717} \tau^{(x,n)}(t) = \inf \biggl\{s>0\dvt \int
_0^{s} \exp \bigl\{ \alpha\xi_u^{(x,n)}
\bigr\} \,\mathrm{d}u > t \biggr\}.
\end{equation}
To verify this, we take $x>0$ and $n$ even, the other cases can be
proved similarly. In this case, $X_{H_n}^{(x)}>0$. By the strong Markov
property, conditionally on $X_{H_n} = y$, we have
\[
( X_{H_n+t}, 0\leq t< H_{n+1} - H_n ) \overset{
\mathcal{L}}= \bigl\{ ( X_t, 0\leq t<H_1 ),
\mathbb{P}_y \bigr\}.
\]
And since the process on the right-hand side of the latter expression
is a positive self-similar Markov process, then by Lamperti's
representation there exists a Levy process $(\xi^+,\mathbf{P})$ such that
\[
\bigl\{ ( X_t, 0\leq t<H_1 ), \mathbb{P}_y
\bigr\} \overset {\mathcal{L}}= \bigl\{ \bigl( y\exp \bigl\{
\xi_{\tau^+(ty^{-\alpha})}^+ \bigr\}, 0\leq t< A^+(\infty) \bigr), \mathbf{P} \bigr\},
\]
where
\[
A^+(\infty) = \int_0^\infty\exp \bigl\{ \alpha
\xi_s^+ \bigr\}\,\mathrm{d}s.
\]
Furthermore, since $H_1<\infty$, $\mathbb{P}_y$-a.s., then $\xi^+$ is a
killed L\'evy process with lifetime $\zeta^+$, exponentially
distributed with parameter $q^+>0$ and hence
\[
A^+(\infty) = \int_0^{\zeta^+} \exp \bigl\{ \alpha
\xi_s^+ \bigr\}\,\mathrm{d}s.
\]
Note that we chose the superscript $+$ because $\sign(X_{H_n})>0$.

Thus, we have obtained that for all $x>0$, $n$ even,
\[
\bigl\{ ( X_{H_n+t}, 0\leq t<H_{n+1}-H_n ),
\mathbb{P}_x \bigr\} \overset{\mathcal{L}}= \bigl\{ \bigl(
X_{H_n}^{(x)} \exp \bigl\{ \xi_{\tau^+(t(X_{H_n}^{(x)})^{-\alpha})}^+ \bigr\}, 0
\leq t< A^+(\infty) \bigr), \mathbf{P} \bigr\}.
\]
This shows (\ref{eq1616}). Also, the Lamperti representation ensures
that for all $x\in\mathbb{R}^\ast$, $n\geq0$,
%
\begin{equation}
\label{eq1818} \bigl|X_{H_n}^{(x)}\bigr|^{-\alpha}
(H_{n+1}-H_n) = \int_0^{\zeta^{(x,n)}}
\exp \bigl\{ \alpha\xi_u^{(x,n)} \bigr\} \,\mathrm{d}u,
\end{equation}
which implies that for all $n\geq1$
%
\begin{equation}
\label{eq1919} H_n^{(x)} = \sum
_{k=0}^{n-1}\bigl |X_{H_k}^{(x)}\bigr|^{\alpha}
\int_0^{\zeta
^{(x,k)}} \exp \bigl\{ \alpha
\xi_u^{(x,k)} \bigr\} \,\mathrm{d}u.
\end{equation}
Now, for $x \in\mathbb{R}^\ast$ we define the sequence $(U^{(x,n)},
n\geq0)$ by
\[
\exp \bigl\{ U^{(x,n)} \bigr\} = -\frac{ X_{H_{n+1}}^{(x)} }{ X_{H_{n+1}-}^{(x)}
},\qquad  n\geq0.
\]
Then, by (\ref{eq1616}) and (\ref{eq1818}) it follows that
\[
X_{H_{n+1}-}^{(x)} = X_{H_n}^{(x)} \exp \bigl\{
\xi_{\zeta}^{(x,n)} \bigr\},
\]
and also
\[
X_{H_{n+1}}^{(x)} = X_{H_{n+1}-}^{(x)}
\frac{ X_{H_{n+1}}^{(x)} }{
X_{H_{n+1}-}^{(x)} } = - X_{H_n}^{(x)} \exp \bigl\{
\xi_{\zeta}^{(x,n)} + U^{(x,n)} \bigr\}.
\]
Hence, for all $n\geq0$,
%
\begin{equation}
\label{eq2020} X_{H_{n+1}}^{(x)} = x\exp \Biggl\{ \sum
_{k=0}^n \bigl( \xi_{\zeta
}^{(x,k)}
+ U^{(x,k)} \bigr) + \mathrm{i}\ppi(n+1) \Biggr\}.
\end{equation}
Note that because of Theorem \ref{theo4}, for every $x\in\mathbb
{R}^\ast$, the sequence $(\xi^{(x,n)}, \zeta^{(x,n)}, U^{(x,n)}, n\geq
0)$ satisfies the condition which defines the process $Y^{(x)}$ in (\ref
{eq1010}). It only remains to prove that $X^{(x)}$ time changed is of
the form (\ref{eq1010}). For that aim, write
\[
A^{(x,n)}(t) = \int_0^t \exp \bigl\{
\alpha\xi_u^{(x,n)} \bigr\}\,\mathrm{d}u,\qquad  0\leq t\leq\zeta^{(x,n)}.
\]
Thanks to (\ref{eq1616}), (\ref{eq1717}) and (\ref{eq2020}), we have
\begin{eqnarray*}
X_{H_n + |X_{H_n}^{(x)}|^{\alpha} A^{(x,n)}(t)}^{(x)} &=& X_{H_n}^{(x)} \exp \bigl
\{ \xi_t^{(x,n)} \bigr\}
\\
&=& x\exp \bigl\{ \mathcal{E}_{t+T_n}^{(x)} \bigr\}.
\end{eqnarray*}
On the other hand, by (\ref{eq1919}), for $0\leq t< \zeta^{(x,n)}$ it follows
\begin{eqnarray*}
&&H_n^{(x)} + \bigl|X_{H_n}^{(x)}\bigr|^{\alpha}
A^{(x,n)}(t)
\\
&&\quad= \sum_{k=0}^{n-1} \bigl|X_{H_k}^{(x)}\bigr|^{\alpha}
\int_0^{\zeta^{(x,k)}} \exp \bigl\{ \alpha
\xi_u^{(x,k)} \bigr\} \,\mathrm{d}u + \bigl|X_{H_n}^{(x)}\bigr|^{\alpha}
\int_0^{t} \exp \bigl\{ \alpha
\xi_u^{(x,n)} \bigr\} \,\mathrm{d}u
\\
&&\quad= \sum_{k=0}^{n-1} \int
_0^{\zeta^{(x,k)}} |x|^\alpha\bigl|\exp \bigl\{ \alpha
\mathcal{E}_{u+T_k}^{(x)} \bigr\}\bigl| \,\mathrm{d}u + \int_0^{t}
|x|^\alpha\bigl|\exp \bigl\{ \alpha\mathcal{E}_{u+T_n}^{(x)}
\bigr\}\bigl| \,\mathrm{d}u
\\
&&\quad= \sum_{k=0}^{n-1} \int
_{T_k}^{T_{k+1}} |x|^\alpha\bigl|\exp \bigl\{ \alpha
\mathcal{E}_{u}^{(x)} \bigr\}\bigr| \,\mathrm{d}u + \int_{T_n}^{t+T_n}
|x|^\alpha\bigl|\exp \bigl\{ \alpha\mathcal{E}_{u}^{(x)}
\bigr\}\bigr| \,\mathrm{d}u
\\
&&\quad= \int_0^{t+T_n}\bigl |x \exp \bigl\{
\mathcal{E}_{u}^{(x)} \bigr\}\bigr|^\alpha \,\mathrm{d}u.
\end{eqnarray*}
Hence,
\[
X_{\int_0^{t+T_n} |x\exp\{\mathcal{E}_s^{(x)}\}|^\alpha \,\mathrm{d}s}^{(x)} = x\exp \bigl\{\mathcal{E}_{t+T_n}^{(x)}
\bigr\},\qquad 0\leq t<\zeta^{(x,n)}.
\]
The latter and (\ref{eq1515}) imply that $\mathcal{Y}^{(x)}$ can be
decomposed as in (\ref{eq1010}). Furthermore, as a consequence of this
decomposition and the converse of the main result in Kiu \cite{Kiu80}, we
can conclude that every Lamperti--Kiu process can be constructed as
explained in (\ref{eq1010}).

{(iii)} Let $(\mathcal{G}_t)$ be the natural filtration of
$Y^{(x)}$, that is, $\mathcal{G}_t = \sigma(Y_s^{(x)}, s\leq t)$, $t\geq
0$. Let $\mathcal{F}_t = \mathcal{G}_{\tau(t|x|^{-\alpha})}$, $t\geq
0$. Clearly, $X^{(x)}$ is $(\mathcal{F}_t)$-adapted, and since the
strong Markov property\vspace*{-1pt} is preserved under time changes by additive
functionals, $X^{(x)}$ is a strong Markov process. We recall $\mathcal
{E}^{(cx)} \overset{\mathcal{L}}= \mathcal{E}^{(x)}$ for all $c>0$.
Thus, if $c>0$, then
\begin{eqnarray*}
\bigl( cX_{c^{-\alpha}t}^{(x)}, t\geq0 \bigr) &=& \bigl( cx\exp \bigl
\{ \mathcal{E}_{\tau(t|cx|^{-\alpha})}^{(x)} \bigr\}, t\geq0 \bigr)
\\
&\overset{\mathcal{L}}=& \bigl( cx\exp \bigl\{ \mathcal{E}_{\tau
(t|cx|^{-\alpha})}^{(cx)}
\bigr\}, t\geq0 \bigr)
\\
&=& \bigl( X_{t}^{(cx)}, t\geq0 \bigr).
\end{eqnarray*}
This proves the self-similar property of $X^{(x)}$. It only remains to
prove that all $X^{(x)}$ have the same semigroup. We have
\[
X_{t+s}^{(x)} = X_t^{(x)}
\bigl(Y_{\tau(t|x|^{-\alpha})}^{(x)} \bigr)^{-1} Y_{\tau
((t+s)|x|^{-\alpha})}^{(x)}.
\]
On the other hand, for all $s,t\geq0$,
\begin{eqnarray*}
&&\tau \bigl((t+s)|x|^{-\alpha} \bigr)
\\
&&\quad= \tau \bigl(t|x|^{-\alpha} \bigr) + \inf \biggl\{r>0\dvt \int
_0^r \bigl|\exp \bigl\{ \alpha \mathcal{E}_{\tau(t|x|^{-\alpha})+u}^{(x)}
\bigr\}\bigr| \,\mathrm{d}u > s|x|^{-\alpha} \biggr\}
\\
&&\quad= \tau \bigl(t|x|^{-\alpha} \bigr) + \inf \biggl\{r>0\dvt \int
_0^r \bigl\llvert \bigl(Y_{\tau
(t|x|^{-\alpha})}^{(x)}
\bigr)^{-1} Y_{\tau(t|x|^{-\alpha})+u}^{(x)} \bigr\rrvert^{\alpha}
\,\mathrm{d}u > s\bigl|X_t^{(x)}\bigr|^{-\alpha} \biggr\}.
\end{eqnarray*}
Write $\widehat{Y}_s^{(x)} = (Y_{\tau(t|x|^{-\alpha})}^{(x)})^{-1}
Y_{\tau(t|x|^{-\alpha})+s}^{(x)}$, $s\geq0$. Then
\[
X_{t+s}^{(x)} = X_t^{(x)}
\widehat{Y}_{\widehat{\tau
}(s|X_t^{(x)}|^{-\alpha})}^{(x)}.
\]
Hence by the strong Markov property of $Y^{(x)}$, Theorem \ref{theo7}{(ii)}, we obtain
\begin{eqnarray*}
\mathbf{P} \bigl( X_{t+s}^{(x)} \in \mathrm{d}z \mid
\mathcal{F}_t \bigr) &=& \mathbf{P} \bigl( X_t^{(x)}
\widehat{Y}_{ \widehat{\tau}(s|X_t^{(x)}|^{-\alpha})
}^{(x)}\in \mathrm{d}z \mid\mathcal{F}_t
\bigr)
\\
&=& \mathbf{P} \bigl( y \exp \bigl\{\mathcal{E}_{ \tau(s|y|^{-\alpha}) }^{(\sign
(y))}
\bigr\}\in \mathrm{d}z \bigr)|_{{y=X_t^{(x)}}}
\\
&=& \mathbf{P} \bigl( X_t^{(y)}\in \mathrm{d}z \bigr)|_{{y=X_t^{(x)}}}.
\end{eqnarray*}
This concludes the proof.
\end{pf*}

\begin{remark}
Let $A^{(x)} = (A_t^{(x)}, 0\leq t\leq\infty)$ be the process defined by
\[
A_t^{(x)} = \int_0^t \bigl|
\exp \bigl\{ \alpha\mathcal{E}_s^{(x)} \bigr\}\bigr|\,\mathrm{d}s,\qquad  0\leq t
\leq \infty.
\]
Note that $A^{(x)}$ only depends on $x$ through its sign. From (\ref
{eq1919}), (\ref{eq2020}) and Proposition \ref{prop2}, under $\mathbb{P}_x$,
\[
T = \lim_{n\rightarrow\infty} H_n = |x|^{\alpha}A_\infty^{(\sign(x))},
\]
that is, there is a relation between the hitting time of zero for $X$ and
the exponential functional of~$\mathcal{E}$, similar to the one known
for positive self-similar Markov processes. Furthermore, Lamperti's
representation can be written as
\[
X_t^{(x)} \mathbf{1}_{\{t<T\}} = x\exp \bigl\{
\mathcal{E}_{\tau
^{(x)}(t|x|^{-\alpha})}^{(x)} \bigr\} \mathbf{1}_{\{t<|x|^\alpha A_{\infty
}^{(\sign(x))}\}},\qquad  t
\geq0,
\]
where $\tau^{(x)}(t) = \inf\{s>0\dvt  \int_0^s |\exp\{ \alpha\mathcal
{E}_u^{(x)} \}| \,\mathrm{d}u > t \}$, $t < A_\infty^{(x)}$.
\end{remark}

\begin{pf*}{Proof of Proposition \ref{prop8}}
We prove the case $x>0$, the case $x<0$ can be proved similarly. Let~$T_1$ and
$T_2$ the first and the second times of sign change for $Y$,
respectively. In the case $x>0$,
\[
T_1=\inf\{t>0\dvt Y_t<0\}, \qquad T_2=\inf
\{t>T_1\dvt Y_t>0\}.
\]
Since $f$ is bounded, we have
\[
\mathbf{E}_x \bigl[f(Y_t) \bigr] -f(x) =
\mathbf{E}_x \bigl[f(Y_t)\mathbf{1}_{\{ T_1>t \}
} -
f(x) \bigr] + \mathbf{E}_x \bigl[f(Y_t)
\mathbf{1}_{\{ T_1\leq t <T_2 \}} \bigr] + \mathbf{E}_x
\bigl[f(Y_t) \mathbf{1}_{\{ T_2\leq t \}} \bigr].
\]
Recall that by construction of $Y$, $(T_1, T_2)$ are such that under
$\mathbf{P}_x$, for $x>0$, they have the same distribution as $(\zeta^+, \zeta^+ + \zeta^-)$, with $\zeta^+$, $\zeta^-$ independent
exponential random variables with parameters $q^+$, $q^-$,
respectively. It is easy to verify that
\[
\mathbf{P}_x(T_2\leq t) = \cases{
 \displaystyle{\frac{ q^-(1-\mathrm{e}^{-q^+t}) - q^+(1-\mathrm{e}^{-q^-t}) }{ q^- -
q^+}}, & \quad $q^+\neq q^-$,
\vspace*{2pt}\cr
 \displaystyle{1-\mathrm{e}^{-q^+t}-q^+t\mathrm{e}^{-q^+t}}, & \quad $q^+=q^-$.}
\]
It follows that $\mathbf{P}_x(T_2\leq t)=\mathrm{o} ( q^+q^-t^2/2  )$
as $t\rightarrow0$. Hence, using again that $f$ is bounded, we obtain
\[
\frac{1}{t}\mathbf{E}_x \bigl[f(Y_t)
\mathbf{1}_{\{ T_2\leq t \}} \bigr] \leq\frac
{1}{t} C\mathbf{P}_x(T_2
\leq t) \rightarrow0, \qquad t\rightarrow0.
\]
Now we write
\[
\frac{1}{t} \bigl( \mathbf{E}_x \bigl[f(Y_t)
\mathbf{1}_{\{ T_1>t \}} \bigr] - f(x) \bigr) = \frac{1}{t} \bigl(
\mathbf{E}_x \bigl[f \bigl( \exp \bigl\{\xi_t^+ \bigr\}
\bigr) \bigr] -f(x) \bigr) \mathrm{e}^{-q^+t} + \frac
{1}{t}f(x)
\bigl(\mathrm{e}^{-q^+t}-1 \bigr),
\]
where $\xi^+$ is a L\'evy process such that $\xi_0^+ = \log(x)$,
$\mathbf{P}_x$-a.s. The last expression implies
\[
\lim_{t\rightarrow0} \frac{1}{t} \bigl( \mathbf{E}_x
\bigl[f(Y_t)\mathbf{1}_{\{
T_1>t \}} \bigr] - f(x) \bigr) =
\mathcal{A}^{+}(f\circ\exp) \bigl(\log(x) \bigr) - q^+f(x).
\]
To conclude, observe the identity
\[
\mathbf{E}_x \bigl[f(Y_t)\mathbf{1}_{\{ T_1\leq t <T_2 \}}
\bigr] = \mathbf{E}_x \bigl[ f \bigl(-\exp \bigl\{
\xi_{t-\zeta^+}^- + \xi_{\zeta^+}^+ + U_1^+ \bigr\} \bigr)
\mid0\leq t-\zeta^+< \zeta^- \bigr] \mathbf{P}_x(T_1
\leq t < T_2),
\]
where $\xi^+$ is as before and $\xi^-$ is a L\'evy process with
lifetime $\zeta^-$ independent of $(\xi^+, \zeta^+, U_1^+)$ and
satisfying $\xi_0^-=0$, $\mathbf{P}_x$-a.s. This together with
\[
\lim_{t\rightarrow0}\frac{1}{t}\mathbf{P}_x(T_1
\leq t < T_2) = \lim_{t\rightarrow0}\frac{1}{t}
\mathbf{P}_x(T_1\leq t) - \lim_{t\rightarrow0}
\frac{1}{t}\mathbf{P}_x(T_2\leq t) = q^+,
\]
and the convergence
\[
\lim_{t\rightarrow0} \frac{1}{t} \mathbf{E}_x \bigl[ f
\bigl(-\exp \bigl\{\xi_{t-\zeta
^+}^- + \xi_{\zeta^+}^+ +
U_1^+ \bigr\} \bigr) \mid0\leq t-\zeta^+<\zeta^- \bigr] = \mathbf{E}
\bigl[f \bigl(-x\exp \bigl\{U^+ \bigr\} \bigr) \bigr],
\]
which holds by the right continuity of $\xi^+$ and $\xi^-$, imply that
\[
\lim_{t\rightarrow0}\frac{1}{t}\mathbf{E}_x
\bigl[f(Y_t)\mathbf{1}_{\{
T_1\leq t <T_2 \}} \bigr] = q^+\mathbf{E} \bigl[f
\bigl(-x\exp \bigl\{U^+ \bigr\} \bigr) \bigr].
\]
This ends the proof.
\end{pf*}

\section{Examples} \label{examples}

The aim of this section is to characterize the law of $(\xi^{\pm}, \zeta^{\pm}, U^{\pm})$
which defines the Lamperti--Kiu processes through two
examples. The first example is the $\alpha$-stable process killed at
the first hitting time of zero, and the second is the $\alpha$-stable
process conditioned to avoid zero in the case $\alpha\in(1,2)$.

We start by reviewing some results in the literature about self-similar
Markov processes. Through this section, $X$ will denote an $\alpha
$-stable process and $T$ its first hitting time of zero ($T=\inf\{t>0\dvt
X_t=0\}$, with $\inf\{\varnothing\}=\infty$); and we will denote by $X^0$
and $X^\updownarrow$ the $\alpha$-stable process killed at $T$ and
conditioned to avoid zero, respectively.

In the case $\alpha=2$, the process $X$ has no jumps and $X^0$
corresponds to a standard real Brownian motion absorbed at level 0. On
the other hand, the Brownian motion conditioned to avoid zero is a
three dimensional Bessel process, see, for example, Revuz and Yor \cite{Revuz-Yor}. Thus,
depending on the starting point, $X^\updownarrow$ is such that
$X^\updownarrow$ or $-X^\updownarrow$ is a Bessel process of dimension
3. Since all Bessel processes are obtained as the images by the
Lamperti representation of the exponential of Brownian motion with
drift, see, for example, Carmona, Petit and Yor \cite{Carmona-Petit-Yor} or Yor \cite{Yor92}, we obtain
the following for $x\in\mathbb{R}^\ast$,
\[
X_t^0 = x\exp \bigl\{\xi_{\tau(t|x|^{-\alpha})}^0
\bigr\}, \qquad X_t^\updownarrow= x\exp \bigl\{\xi_{\tau(t|x|^{-\alpha})}^\updownarrow
\bigr\},\qquad  t\geq0,
\]
where $\xi^0$ and $\xi^\updownarrow$ are real Brownian motions with
drift, viz., $\xi^0 = (B_t - t/2, t\geq0)$ and $\xi^\updownarrow=
(\widetilde{B}_t + t/2, t\geq0)$, with $B$, $\widetilde{B}$ real
Brownian motions. Therefore, the Lamperti representation is known in
the case $\alpha=2$, so we exclude this case in our examples.

For $0<\alpha<2$, let $\psi$ be the characteristic exponent of $X$:
$\mathbb{E}[\exp{(\mathrm{i}\lambda X_t)}]=\exp{(t\psi(\lambda))}$, $t\geq0$,
$\lambda\in\mathbb{R}$. It is well known that $\psi$ is given by
%
\begin{equation}
\label{eqaux1} \psi(\lambda) = \mathrm{i}a\lambda+ \int_{\mathbb{R}}
\bigl(\mathrm{e}^{\mathrm{i}\lambda y}-1-\mathrm{i}\lambda y\mathbf{1}_{\{|y|<1\}} \bigr) \nu(y)\,\mathrm{d}y,\qquad
\lambda\in\mathbb{R},
\end{equation}
where $\nu$ is the density of the L\'evy measure:
%
\begin{equation}
\label{eqaux2} \nu(y) = c^+y^{-\alpha-1}\mathbf{1}_{\{y>0\}} +
c^-|y|^{-\alpha
-1}\mathbf{1}_{\{y<0\}},
\end{equation}
with $c^+$ and $c^-$ being two non-negative constants such that
$c^++c^->0$. The constant $a$ is $(c^+-c^-)/(1-\alpha)$ if $\alpha\neq
1$. For the case $\alpha=1,$ we will assume that $X$ is a symmetric
Cauchy process, thus $c^+=c^-$ and $a=0$.

Another quite well studied positive self-similar Markov process killed
at its first hitting time of $0$ is the process obtained by killing an
$\alpha$-stable process when it leaves the positive half-line.
Formally, if $R$ is the stopping time $R=\inf\{t>0\dvt  X_t\leq0\}$, then
the process killed at the first time it leaves the positive half-line
is $X^\dagger= (X_t\mathbf{1}_{\{t<R\}}, t\geq0)$ where $0$ is
assumed to be a cemetery state. Caballero and Chaumont in \cite{Caballero-Chaumont}
proved that the L\'evy process $\xi$ related to
$X$ via Lamperti's representation has the characteristic exponent:
%
\begin{equation}
\label{eq2121} \Phi(\lambda) = \mathrm{i}a\lambda+ \int_{\mathbb{R}}
\bigl[\mathrm{e}^{\mathrm{i}\lambda y} - 1 - \mathrm{i}\lambda \bigl(\mathrm{e}^y-1 \bigr)
\mathbf{1}_{\{|\mathrm{e}^y-1|<1\}} \bigr] \pi(\mathrm{d}y) - c^-\alpha^{-1}, \qquad \lambda\in
\mathbb{R},
\end{equation}
where the L\'evy measure $\pi(\mathrm{d}y)$ is
%
\begin{equation}
\label{eq2222} \pi(\mathrm{d}y) = \biggl( \frac{ c^+ \mathrm{e}^y }{ (\mathrm{e}^y-1)^{\alpha+1} } \mathbf{1}_{\{
y>0 \}} +
\frac{ c^- \mathrm{e}^y }{ (1-\mathrm{e}^y)^{\alpha+1} } \mathbf{1}_{\{ y<0 \}} \biggr) \,\mathrm{d}y.
\end{equation}
Note from (\ref{eq2121}) that the killing rate of the L\'evy process
$\xi$ is $c^-\alpha^{-1}$.

A further example in the literature appears in Caballero, Pardo and P\'erez \cite
{Caballero-Pardo-Perez11}. They studied the radial part of the
symmetric $\alpha$-stable process taking values in $\mathbb{R}^d$. In
the case $d=1$, $0<\alpha<1$, they proved that the L\'evy process in
the Lamperti representation for the radial part of the symmetric $\alpha
$-stable process is the sum of two independent L\'evy processes $\xi_1$, $\xi_2$ with triples $(0,0,\pi_1)$ and $(0,0,\pi_2)$ where
%
\begin{eqnarray}
\label{eq2323} \pi_1(\mathrm{d}y) &=& \biggl( \frac{ k(\alpha) \mathrm{e}^y }{ (\mathrm{e}^y-1)^{\alpha+1} } \mathbf
{1}_{\{ y>0 \}} + \frac{ k(\alpha) \mathrm{e}^y }{ (1-\mathrm{e}^y)^{\alpha+1} } \mathbf {1}_{\{ y<0 \}} \biggr)\,\mathrm{d}y,
\nonumber
\\[-8pt]
\\[-8pt]
\nonumber
\pi_2(\mathrm{d}y) &=& \frac{ k(\alpha) \mathrm{e}^y}{(\mathrm{e}^y+1)^{\alpha+1}} \,\mathrm{d}y
\end{eqnarray}
and
\[
k(\alpha) = \frac{\alpha}{2\Gamma(1-\alpha)\cos{\ppi\alpha}/{2}}.
\]
In other words, the L\'evy process in the Lamperti representation is
the sum of a L\'evy process with L\'evy measure similar to (\ref
{eq2222}) and a compound Poisson process. Since the process $Y$ is
symmetric in this case, the results in Caballero, Pardo and P\'erez
\cite{Caballero-Pardo-Perez11}
confirm Chybiryakov's results.

The L\'evy processes with L\'evy measure having the form (\ref{eq2222})
or $\pi_1$ in (\ref{eq2323}) are examples of Lamperti-stable processes.
For the definition and properties of Lamperti-stable processes, see
Caballero, Pardo and P\'erez \cite{Caballero-Pardo-Perez10}.

\subsection{\texorpdfstring{The $\alpha$-stable process killed at zero}{The alpha-stable process killed at zero}} \label{alpha-stablekilledatzero}

The following theorem provides the expression of the infinitesimal
generator of the process~$X^0$.

\begin{theorem} \label{infgen1}
Let $\alpha\in(0,2)$ and let $\mathcal{A}$, $\mathcal{A}^0$ the
infinitesimal generators of the $\alpha$-stable process and the $\alpha
$-stable process killed in $T$, respectively. Then $\mathcal
{D}_{\mathcal{A}^0} = \{f\in\mathcal{D}_{\mathcal{A}}\dvt  f(0) = 0\}$ and
$\mathcal{A}^0 f(x) = \mathcal{A}f(x)$, for $x\in\mathbb{R}^\ast$.
Furthermore, for $x\in\mathbb{R}^\ast$, $\mathcal{A}^0 f(x)$ can be
written as
%
\begin{eqnarray}
\label{eq2727} \mathcal{A}^0 f(x) &=& \frac{1}{|x|^{\alpha}} \biggl[
\sign(x)axf'(x) + c^{-\sign(x)}\alpha^{-1} \int
_{\mathbb{R}^-} \bigl[f(xu) - f(x) \bigr] g^{0}(u) \,\mathrm{d}u
\nonumber
\\[-8pt]
\\[-8pt]
\nonumber
&&\hspace*{22pt}{} + \int_{\mathbb{R}^+} \bigl[f(xu) - f(x) - xf'(x)
(u-1)\mathbf {1}_{\{|u-1|<1\}} \bigr]\nu^{0,\sign(x)}(u)\,\mathrm{d}u \biggr],
\end{eqnarray}
where
\[
\nu^{0,\sign(x)}(u) = \nu \bigl(\sign(x) (u-1) \bigr),\qquad u>0,\qquad
g^0(u)= \alpha (1-u)^{-\alpha-1}, \qquad u<0,
\]
and $\nu$ is given by (\ref{eqaux2}).
\end{theorem}

The proof of the latter theorem will be given at the end of this
subsection. The following corollary characterizes the Lamperti--Kiu
process associated to the $\alpha$-stable process killed at its first
hitting time of zero and its proof is an immediate consequence of
Volkonskii's theorem and the formulas (\ref{eq1212}) and (\ref{eq2727}).

\begin{corollary} \label{coro16}
Let $\xi^{0,\pm}, \zeta^{0,\pm}, U^{0,\pm}$ the random objects in the
Lamperti representation of $X^0$. Then, the characteristic exponent of
$\xi^{0,\pm}$ is given by
\[
\psi^{0,\pm}(\lambda) = \mathrm{i}a^{\pm}\lambda+ \int
_{\mathbb{R}} \bigl[\mathrm{e}^{\mathrm{i}\lambda y} - 1 - \mathrm{i}\lambda
\bigl(\mathrm{e}^{y}-1 \bigr)\mathbf{1}_{\{ |\mathrm{e}^y-1|<1\}} \bigr]
\pi^{0,\pm}(\mathrm{d}y), \qquad \lambda\in\mathbb{R},
\]
where $a^{\pm} = \pm a$, with $a$ as in (\ref{eqaux1}), and $\pi^{0,\pm
}(\mathrm{d}y) = \mathrm{e}^y\nu(\pm(\mathrm{e}^y-1))\,\mathrm{d}y$. The parameters of the exponential random
variables $\zeta^{0,\pm}$ are $c^{\mp}\alpha^{-1}$ and the real random
variables $U^{0,\pm}$ have density
\[
g(u) = \frac{\alpha \mathrm{e}^u}{(1+\mathrm{e}^u)^{\alpha+1}},\qquad  u\in\mathbb{R}.
\]
\end{corollary}

Note that as expected, the L\'evy process $\xi^{0,+}$ is the one
obtained in Caballero  and Chaumont~\cite{Caballero-Chaumont}. Furthermore, the downwards
change of sign rate, which is the death rate in Caballero  and Chaumont \cite
{Caballero-Chaumont}, is $c^{-}\alpha^{-1}$. From the triples of $\xi^{0,+}$ and $\xi^{0,-}$, we can observe that both belong to the
Lamperti-stable family. In the particular case where $X$ is a symmetric
$\alpha$-stable process with $\alpha\in(0,1)$, the description in
Corollary \ref{coro16} coincides with the one in Caballero, Pardo  and P\'erez \cite
{Caballero-Pardo-Perez11}, see (\ref{eq2323}). Note that $U^{0,+}$,
$U^{0,-}$ are identically distributed and they are such that $U^{0,\pm}
\overset{\mathcal{L}}= \log V$, where $V$ follows a Pareto distribution
with parameter $\alpha$, viz.,
\[
f(x) = \frac{\alpha}{(1+x)^{\alpha+1}},\qquad  x>0.
\]

In order to prove the main theorem of this subsection, we need the
following two lemmas.

\begin{lemma}\label{lemma13}
Let $X$ be an $\alpha$-stable process, $\alpha\in(0,2)$. Then, for any
$x\in\mathbb{R}^\ast$,
%
\begin{equation}
\label{eq2424} \lim_{t\downarrow0} \frac{1}{t} \mathbb{P}_x
\bigl(T\leq t, X_t\in\mathbb {R}^\ast \bigr) = 0.
\end{equation}
\end{lemma}
\begin{pf}
Since for $\alpha\in(0,1]$ the point zero is polar, then (\ref{eq2424})
is clearly satisfied. Suppose $\alpha\in(1,2)$. For $\delta>0$, write
\[
\mathbb{P}_x \bigl(T\leq t, X_t\in
\mathbb{R}^\ast \bigr) = \mathbb{P}_x\bigl(T\leq t,
|X_t|\in(0,\delta]\bigr) + \mathbb{P}_x\bigl(T\leq t,
|X_t|>\delta\bigr).
\]
First, we verify the following: for $0<\delta<|x|$ it holds
%
\begin{equation}
\label{eq2525} \lim_{t\downarrow0} \frac{1}{t} \mathbb{P}_x\bigl(|X_t|
\in(0,\delta]\bigr) = \frac{ c^{-\sign(x)} }{ \alpha} \sign(x) \bigl(|\delta-x|^{-\alpha} - |
\delta +x|^{-\alpha} \bigr).
\end{equation}
For this aim, we will use the fact that for every $K>0$, $(1/t)\mathbb
{P}_0(X_t\in \mathrm{d}z)$ converges vaguely to $\nu(z)\,\mathrm{d}z$ on $\{z\dvt |z|>K\}$, as
$t\downarrow0$; see, for example, exercise I.1 in Bertoin \cite{Bertoin}. We only
show (\ref{eq2525}) in the case $x<0$, the case $x>0$ can be proved
similarly. For $x<0$, we have $\delta+x<0$ and
\begin{eqnarray*}
\lim_{t\downarrow0} \frac{1}{t} \mathbb{P}_x\bigl(|X_t|
\in(0,\delta]\bigr) &=& \lim_{t\downarrow0} \frac{1}{t} \mathbb{P}_0
\bigl(X_t \in[-\delta -x,\delta-x] \bigr)
\\
&=& \int_{-\delta-x}^{\delta-x} \nu(z)\,\mathrm{d}z
\\
&=& \frac{c^+}{\alpha} \bigl((-\delta-x)^{-\alpha} - (\delta-x)^{-\alpha}
\bigr),
\end{eqnarray*}
which proves the claim. Now, from (\ref{eq2525}), we obtain
%
\begin{equation}
\label{eq2626} \limsup_{t\downarrow0} \frac{1}{t} \mathbb{P}_x\bigl(T
\leq t, |X_t|\in (0,\delta]\bigr) \leq\frac{ c^{-\sign(x)} }{ \alpha} \sign(x) \bigl(|
\delta -x|^{-\alpha} - |\delta+x|^{-\alpha} \bigr).
\end{equation}
On the other hand, by the strong Markov property
\[
\mathbb{P}_x\bigl(T\leq t, |X_t|>\delta\bigr) = \int
_0^t \mathbb {P}_0\bigl(|X_{t-s}|>
\delta\bigr) \mathbb{P}_x(T\in \mathrm{d}s).
\]
Since $(1/t)\mathbb{P}_0(X_t\in \mathrm{d}z)$ converges vaguely to $\nu(z)\,\mathrm{d}z$ on
$\{z\dvt  |z|>K\}$ for every $K>0$, there exists a constant $C$ such that,
for sufficiently small $t$:
\[
\mathbb{P}_0\bigl(|X_{t-s}|>\delta\bigr) \leq\frac{Ct}{\delta^\alpha},\qquad
\mbox{for all }s\in(0,t).
\]
Then
\[
\mathbb{P}_x\bigl(T\leq t, |X_t|>\delta\bigr) \leq
\mathbb{P}_x(T\leq t) \frac
{Ct}{\delta^\alpha}.
\]
The latter inequality and (\ref{eq2626}) imply the result.
\end{pf}

\begin{lemma} \label{lemma14}
Let $x\in\mathbb{R}^\ast$, and $\alpha\in(0,2)$. We will denote by
$I_1^{(x)}$ and $I_2^{(x)}$ the following integrals
\begin{eqnarray*}
I_1^{(x)} &=& \int_{\mathbb{R}^+}(u-1) (
\mathbf{1}_{\{|u-1|<1\}} - \mathbf{1}_{\{|x(u-1)|<1\}}) \nu \bigl(\sign(x) (u-1)
\bigr)\,\mathrm{d}u,
\\
I_2^{(x)} &=& \int_{\mathbb{R}^-}(u-1)
\mathbf{1}_{\{|x(u-1)|<1\}} \nu \bigl(\sign(x) (u-1) \bigr)\,\mathrm{d}u.
\end{eqnarray*}
The identity
\[
I_1^{(x)} - I_2^{(x)} = \sign(x)a
\bigl(1-|x|^{\alpha-1} \bigr),\qquad \mbox{holds}.
\]
\end{lemma}
\begin{pf}
We will show the case $x<0$, and $\alpha\neq1$, the other cases can be
proved similarly. First, observe that $|u-1|<1$ if only if $0<u<2$.
Thus, if $x=-1$, then $I_1^{(x)}=I_2^{(x)}=0$ and the lemma is
satisfied. Now, suppose that $-1<x<0$, then $1+x^{-1}<0<2<1-x^{-1}$,
\[
I_1^{(x)} = -\int_{2}^{1-x^{-1}}
c^-(u-1)^{-\alpha} \,\mathrm{d}u = \frac
{c^-}{1-\alpha} \bigl[1-(-x)^{\alpha-1}
\bigr]
\]
and
\[
I_2^{(x)} = -\int_{1+x^{-1}}^{0}
c^+(1-u)^{-\alpha} \,\mathrm{d}u = \frac
{c^+}{1-\alpha} \bigl[1-(-x)^{\alpha-1}
\bigr].
\]
Hence, $I_1^{(x)}-I_2^{(x)} = -a[1-(-x)^{\alpha-1}]$. Finally, suppose
that $x<-1$. In this case, we have $0<1+x^{-1}<1<1-x^{-1}<2$,
$I_2^{(x)}=0$ and
\[
I_1^{(x)} = -\int_{0}^{1+x^{-1}}
c^+(1-u)^{-\alpha} \,\mathrm{d}u + \int_{1-x^{-1}}^{2}
c^-(u-1)^{-\alpha} \,\mathrm{d}u = -a \bigl[1-(-x)^{\alpha-1} \bigr].
\]
This ends the proof.
\end{pf}

\begin{pf*}{Proof of Theorem \ref{infgen1}}
For any $f$ bounded function such that $f(0)=0$, we have for $x\in
\mathbb{R}^\ast$
\[
\mathbb{E}_x \bigl[f \bigl(X_t^0 \bigr) -
f(x) \bigr] = \mathbb{E}_x \bigl[f(X_t) - f(x) \bigr] -
\mathbb {E}_x \bigl[ f(X_t)\mathbf{1}_{\{ T\leq t\}}
\bigr].
\]
On the other hand, by the Lemma \ref{lemma13},
\[
\lim_{t\rightarrow0} \frac{1}{t} \mathbb{E}_x
\bigl[f(X_t)\mathbf{1}_{\{
T\leq t\}} \bigr] = 0.
\]
Then
\[
\lim_{t\rightarrow0} \frac{1}{t} \mathbb{E}_x \bigl[f
\bigl(X_t^0 \bigr) - f(x) \bigr] = \lim_{t\rightarrow0}
\frac{1}{t} \mathbb{E}_x \bigl[f(X_t) - f(x)
\bigr].
\]
Hence, $\mathcal{D}_{\mathcal{A}^0} = \{f\in\mathcal{D}_{\mathcal{A}}\dvt
f(0) = 0\}$ and $\mathcal{A}^0 f(x) = \mathcal{A}f(x)$.

Now we will obtain (\ref{eq2727}). By the first part of the theorem, we
have that for $x\in\mathbb{R}^\ast$, $\mathcal{A}^0 f(x)$ is given by
%
\begin{equation}
\label{eq2828} \mathcal{A}^0 f(x) = af'(x) + \int
_{\mathbb
{R}} \bigl[f(x+y)-f(x)-yf'(x)
\mathbf{1}_{\{|y|<1\}} \bigr] \nu(y)\,\mathrm{d}y.
\end{equation}
Let $I$ be the integral in (\ref{eq2828}). Then, with the change of
variables $y=x(u-1)$ we obtain
\begin{eqnarray*}
I &= &\frac{1}{|x|^{\alpha}} \int_{\mathbb{R}} \bigl[f(xu) - f(x) -
xf'(x) (u-1)\mathbf{1}_{\{|x(u-1)|<1\}} \bigr]\nu \bigl(\sign(x)
(u-1) \bigr)\,\mathrm{d}u
\\
&=& \frac{1}{|x|^{\alpha}} \biggl[ \int_{\mathbb{R}^+} \bigl[f(xu) - f(x)
- xf'(x) (u-1)\mathbf{1}_{\{|u-1|<1\}} \bigr]\nu \bigl(\sign(x)
(u-1) \bigr)\,\mathrm{d}u
\\
& &\hspace*{22pt}{}+ \int_{\mathbb{R}^+} \bigl[xf'(x) (u-1) (
\mathbf{1}_{\{|u-1|<1\}} - \mathbf{1}_{\{|x(u-1)|<1\}}) \bigr]\nu \bigl(\sign(x)
(u-1) \bigr)\,\mathrm{d}u
\\
&&\hspace*{22pt}{} + \int_{\mathbb{R}^-} \bigl[f(xu) - f(x)- xf'(x)
(u-1)\mathbf {1}_{\{|x(u-1)|<1\}} \bigr]\nu \bigl(\sign(x) (u-1) \bigr)\,\mathrm{d}u \biggr].
\\
\end{eqnarray*}
With the help of Lemma \ref{lemma14}, we can write $I$ as follows
\begin{eqnarray*}
I &=& \frac{1}{|x|^{\alpha}} \biggl[ \sign(x)axf'(x) + \int
_{\mathbb
{R}^-} \bigl[f(xu) - f(x) \bigr]\nu^{0,\sign(x)}(u)\,\mathrm{d}u -
a|x|^{\alpha}f'(x)
\\
&&\hspace*{22pt}{} + \int_{\mathbb{R}^+} \bigl[f(xu) - f(x) - xf'(x)
(u-1)\mathbf {1}_{\{|u-1|<1\}} \bigr]\nu^{0,\sign(x)}(u)\,\mathrm{d}u \biggr] .
\end{eqnarray*}
Hence, we have
\begin{eqnarray*}
\mathcal{A}^0 f(x) &=& \frac{1}{|x|^{\alpha}} \biggl[ \sign(x)
axf'(x) + \int_{\mathbb{R}^-} \bigl[f(xu) - f(x) \bigr]
\nu^{0,\sign(x)}(u) \,\mathrm{d}u
\\
&&\hspace*{22pt}{} + \int_{\mathbb{R}^+} \bigl[f(xu) - f(x) - xf'(x)
(u-1)\mathbf {1}_{\{|u-1|<1\}} \bigr]\nu^{0,\sign(x)}(u) \,\mathrm{d}u \biggr].
\end{eqnarray*}
Finally, note that\vspace*{-6pt}
\[
\frac{\nu^{0,\sign(x)}(u)}{c^{-\sign(x)}\alpha^{-1}} = g^0(u),\qquad  u<0.
\]
This ends the proof.
\end{pf*}

\subsection{\texorpdfstring{The $\alpha$-stable process conditioned to avoid zero}
{The alpha-stable process conditioned to avoid zero}}

In Yano \cite{Yano10} symmetric L\'evy processes conditioned to avoid zero
were studied. One of the main results in Yano \cite{Yano10} can be stated
as follows. Let $X$ be a L\'evy process with characteristic exponent
$\psi$. Consider the following assumptions
\begin{longlist}[{H.1}]
\item[{H.1}] The origin is regular for itself and $X$ is not a
compound Poisson process.

\item[{H.2}] $X$ is symmetric.
\end{longlist}
Then, under {H.1} and {H.2} the function $h$, given by
\[
h(x) = \frac{1}{\ppi} \int_0^\infty
\frac{1-\cos\lambda x}{\theta
(\lambda)} \,\mathrm{d}\lambda,\qquad  x\in\mathbb{R},
\]
where $\theta(\lambda)=-\operatorname{Re}(\psi(\lambda))$, is an invariant function
with respect to the semigroup, $P_t^0$, of the process $X$ killed at
$T$, the first hitting time of 0. Note that if $X$ is an $\alpha
$-stable process with $\alpha\in(0,2)$, {H.1} and {H.2}
are satisfied if and only if $X$ is symmetric and $\alpha\in(1,2)$. In
this case, the characteristic exponent is given by $\psi(\lambda) =
-|\lambda|^{\alpha}$, $h$ has an explicit form, namely
\[
h(x) = C(\alpha) |x|^{\alpha-1},\qquad  x\in\mathbb{R},
\]
where
\[
C(\alpha)=\frac{\Gamma(2-\alpha)}{\uppi(\alpha-1)} \sin\frac{\alpha\ppi}{2}.
\]
In Pant\'i \cite{Panti12} a generalization of the latter fact is considered.
There it is proved that for $X$ $\alpha$-stable process with $1<\alpha
<2$, the function $h$ given by
%
\begin{equation}
\label{eq2929} h(x) = K(\alpha) \bigl(1-\beta\sign(x) \bigr)|x|^{\alpha-1},\qquad
x \in\mathbb{R},
\end{equation}
where
\[
K(\alpha) = \frac{\Gamma(2-\alpha)\sin(\alpha\ppi/2)}{c\ppi(\alpha
-1)(1+\beta^2\tan^2(\alpha\ppi/2))},
\]
and
%
\begin{equation}
\label{eqaux3} c = -\frac{(c^++c^-)\Gamma(2-\alpha)}{\alpha(\alpha-1)}\cos(\alpha\ppi /2),\qquad  \beta=
\frac{c^+-c^-}{c^++c^-},
\end{equation}
is an invariant function for the semigroup of $X^0$. In fact, this
result is a consequence of a more general result that has been proved
in Pant\'i \cite{Panti12} under the sole assumption {H.1}. Since $h$ is
invariant for the semigroup $P_t^0$ and $h(x)\neq0$, for $x\in\mathbb
{R}^\ast$, then we define the semigroup $P_t^h$ on $\mathbb{R}^\ast$ by
\[
P_t^h(x, \mathrm{d}y) := \frac{h(y)}{h(x)}P_t^0(x,\,\mathrm{d}y),\qquad
x,  y\in\mathbb {R}^\ast, t\geq0.
\]
We denote by $\mathbb{P}_x^h$ the law of the strong Markov process with
starting point $x$ and semigroup $P_t^h$. $\mathbb{P}_{\cdot}^h$ is
Doob's h-transformation of $\mathbb{P}_0$ via the invariant function
$h$ as defined in (\ref{eq2929}). Since under the measure $\mathbb
{P}_x^h$ it holds $\mathbb{P}_x^h(T=\infty)=1$, then the process $X^h$
can be considered as the process $X$ conditioned to avoid (or never to
hit) zero, this has been proved in Pant\'i \cite{Panti12}. We use the notation
$X^\updownarrow$ instead of $X^h$ to emphasize this fact. Thus, as was
mentioned at the beginning of the section, $X^\updownarrow$ is the
$\alpha$-stable process conditioned to avoid zero, when $\alpha\in
(1,2)$. In the following lemma, we summarize properties of the function
$h$, which follow straightforwardly from its definition and so we omit
their proof.

\begin{lemma} \label{lemma17}
The function $h$ defined in (\ref{eq2929}) satisfies the following properties
\begin{enumerate}[(iii)]
\item[(i)] $h(x)>0$, for all $x\in\mathbb{R}^\ast$, $h(0)=0$;

\item[(ii)] $h(ux) = |u|^{\alpha-1} h(\sign(u)x)$, for all $u\in\mathbb{R}$;

\item[(iii)] $(hf)'(x) = h(x)[(\alpha-1)x^{-1}f(x)+f'(x)]$, $f\in C^1$,
$x\in\mathbb{R}^\ast$;

\item[(iv)] $h(-x) = h(x) + 2K(\alpha) \beta\sign(x)|x|^{\alpha-1}$,
for all $x\in\mathbb{R}$.
\end{enumerate}
\end{lemma}

Using {(ii)} of Lemma \ref{lemma17} and (\ref{eq111}), it is
possible to verify that the semigroup of the process $X^\updownarrow$
satisfies the self-similarity property. Hence, $X^\updownarrow$ is
real-valued self-similar Markov process. The following theorem provides
an expression for the infinitesimal generator of~$X^\updownarrow$.

\begin{theorem} \label{infgen2}
Let $\mathcal{A}^\updownarrow$ be the infinitesimal generator of
$X^\updownarrow$. For $x\in\mathbb{R}^\ast$, $\mathcal{A}^\updownarrow
f(x)$ can be written as
%
\begin{eqnarray}
\label{eq3030} \mathcal{A}^\updownarrow f(x) &=& \frac{1}{|x|^{\alpha}} \biggl[
a^{\updownarrow, \sign(x)}xf'(x) + c^{\sign(x)}\alpha^{-1}\int
_{\mathbb
{R}^-} \bigl[ f(xu) - f(x) \bigr] g^\updownarrow(u) \,\mathrm{d}u
\nonumber
\\[-8pt]
\\[-8pt]
\nonumber
&&\hspace*{22pt}{} + \int_{\mathbb{R}^+} \bigl[ f(xu) - f(x) - xf'(x)
(u-1)\mathbf {1}_{\{|u-1|<1\}} \bigr]\nu^{\updownarrow, \sign(x)}(u) \,\mathrm{d}u \biggr],
\end{eqnarray}
where
%
\begin{equation}
\label{eq3131} a^{\updownarrow, \sign(x)} = \sign(x)a + c^{\sign(x)}\int
_0^1 \frac{
(1+u)^{\alpha-1}-1 }{ u^\alpha} \,\mathrm{d}u - c^{-\sign(x)}
\int_0^1 \frac{
(1-u)^{\alpha-1} - 1 }{ u^\alpha} \,\mathrm{d}u
\end{equation}
and
\begin{eqnarray*}
\nu^{\updownarrow, \sign(x)}(u)& = &u^{\alpha-1}\nu \bigl(\sign(x) (u-1) \bigr),\qquad
u>0,\\
g^\updownarrow(u)& =& \alpha(-u)^{\alpha-1} (1-u)^{-\alpha
-1},\qquad u<0.
\end{eqnarray*}
\end{theorem}

The following corollary is also a consequence of Volkonskii's theorem
and the comparison of~(\ref{eq1212}) and (\ref{eq3030}).

\begin{corollary}
Let $\xi^{\updownarrow, \pm}, U^{\updownarrow, \pm}, \zeta^{\updownarrow
, \pm}$ the random objects in the Lamperti representation of
$X^\updownarrow$. Then the characteristic exponent of $\xi^{\pm}$ is
\[
\psi^{\updownarrow, \pm}(\lambda) = \mathrm{i}a^{\updownarrow,\pm}\lambda+ \int
_{\mathbb{R}} \bigl[\mathrm{e}^{\mathrm{i}\lambda y} - 1 - \mathrm{i}\lambda
\bigl(\mathrm{e}^{y}-1 \bigr)\mathbf{1}_{\{
|\mathrm{e}^y-1|<1\}} \bigr]
\pi^{\updownarrow, \pm}(\mathrm{d}y),\qquad  \lambda\in\mathbb{R},
\]
where $a^{\updownarrow, \pm}$ is given by (\ref{eq3131}) and $\pi^{\updownarrow, \pm}(\mathrm{d}y) = \mathrm{e}^{\alpha y} \nu(\pm(\mathrm{e}^y-1))\,\mathrm{d}y$. The
parameters of the exponential random variables $\zeta^{\updownarrow, \pm
}$ are $c^{\pm}\alpha^{-1}$ and the real random variables
$U^{\updownarrow, \pm}$ have density
\[
g(u) = \frac{\alpha \mathrm{e}^{\alpha u}}{(1+\mathrm{e}^u)^{\alpha+1}},\qquad u\in\mathbb{R}.
\]
\end{corollary}

As in the first example, the L\'evy processes $\xi^{\updownarrow,+}$,
$\xi^{\updownarrow, -}$ belong to the Lamperti-stable family.
Furthermore, their L\'evy measure, satisfy the relation: $\pi^{\updownarrow, \pm}(\mathrm{d}y) = \mathrm{e}^{(\alpha-1)y}\pi^{0,\pm}(\mathrm{d}y)$. Note that
$g(u)$ can be written as
\[
g(u) = \frac{\alpha \mathrm{e}^{-u}}{(1+\mathrm{e}^{-u})^{\alpha+1}},\qquad  u\in\mathbb{R}.
\]
Hence, $U^{\updownarrow, \pm} \overset{\mathcal{L}}= -U^{0,\pm} \overset
{\mathcal{L}}= -\log V$, with $U^{0,\pm}$ as in Corollary \ref{coro16}
and $V$ is a Pareto random variable.

\begin{pf*}{Proof of Theorem \ref{infgen2}}
Recall that $\mathcal{A}^\updownarrow f(x) = [h(x)]^{-1}\mathcal
{A}^0(hf)(x)$, $x\in\mathbb{R}^\ast$. Thus, by (\ref{eq2727}) we can
write for $x\in\mathbb{R}^\ast$
\[
\bigl[h(x) \bigr]^{-1}|x|^\alpha\mathcal{A}^0(hf)
(x) = \bigl[h(x) \bigr]^{-1} \bigl(\sign (x)ax(hf)'(x) +
\mathcal{I}_1^{(x)} + \mathcal{I}_2^{(x)}
\bigr),
\]
where
\begin{eqnarray*}
\mathcal{I}_1^{(x)} &=&\int_{\mathbb{R}^+}
\bigl[ (hf) (xu) - (hf) (x) - x(hf)'(x) (u-1)\mathbf{1}_{\{|u-1|<1\}}
\bigr] \nu \bigl(\sign(x) (u-1) \bigr)\,\mathrm{d}u,
\\
\mathcal{I}_2^{(x)} &=& \int_{\mathbb{R}^-}
\bigl[ (hf) (xu) -(hf) (x) \bigr] \nu \bigl(\sign(x) (u-1) \bigr)\,\mathrm{d}u.
\end{eqnarray*}
Now, by {(iii)} of Lemma \ref{lemma17},
%
\begin{equation}
\label{eq3232} \bigl[h(x) \bigr]^{-1}\sign(x)ax(hf)'(x)
= \sign(x)axf'(x) + \sign(x)a(\alpha-1)f(x).
\end{equation}
Also, using {(ii)} and {(iii)} of Lemma \ref{lemma17}, we have
%
\begin{eqnarray}
\label{eq3333}
\nonumber
\bigl[h(x) \bigr]^{-1}\mathcal{I}_1^{(x)}
&=& \int_{\mathbb{R}^+} \bigl[ f(xu) - f(x) - xf'(x)
(u-1)\mathbf{1}_{\{|u-1|<1\}} \bigr]u^{\alpha-1}\nu \bigl(\sign(x) (u-1)
\bigr)\,\mathrm{d}u\\
& &{} + \int_{\mathbb{R}^+} \bigl(u^{\alpha-1}-1 \bigr)
(u-1)\mathbf{1}_{\{|u-1|<1\}} \nu \bigl(\sign(x) (u-1) \bigr)\,\mathrm{d}u \times
xf'(x)
\nonumber
\\[-8pt]
\\[-8pt]
\nonumber
& & {}+ \int_{\mathbb{R}^+} \bigl[ u^{\alpha-1}-1-(
\alpha-1) (u-1)\mathbf{1}_{\{
|u-1|<1\}} \bigr] \nu \bigl(\sign(x) (u-1) \bigr)\,\mathrm{d}u
\times f(x)\qquad
\\
&=& I_1^{(x)} + I_2^{(x)}xf'(x)
+ I_3^{(x)}f(x),\nonumber
\end{eqnarray}
where
\begin{eqnarray*}
I_1^{(x)} &=& \int_{\mathbb{R}^+} \bigl[ f(xu)
- f(x) - xf'(x) (u-1)\mathbf {1}_{\{|u-1|<1\}}
\bigr]u^{\alpha-1}\nu \bigl(\sign(x) (u-1) \bigr)\,\mathrm{d}u,
\\
I_2^{(x)} &=& \int_{\mathbb{R}^+}
\bigl(u^{\alpha-1}-1 \bigr) (u-1)\mathbf{1}_{\{
|u-1|<1\}} \nu \bigl(\sign(x)
(u-1) \bigr)\,\mathrm{d}u
\\
&=& c^{\sign(x)}\int_0^1
\frac{ (1+u)^{\alpha-1}-1 }{ u^\alpha} \,\mathrm{d}u - c^{-\sign(x)} \int_0^1
\frac{ (1-u)^{\alpha-1} - 1 }{ u^\alpha} \,\mathrm{d}u,
\\
I_3^{(x)} &=& \int_{\mathbb{R}^+} \bigl[
u^{\alpha-1}-1-(\alpha -1) (u-1)\mathbf{1}_{\{|u-1|<1\}} \bigr] \nu \bigl(
\sign(x) (u-1) \bigr)\,\mathrm{d}u.
\end{eqnarray*}
And by {(ii)}, {(iv)} of Lemma \ref{lemma17} and since $\int_{\mathbb{R}^-} (-u)^{\alpha-1} \nu(\sign(x)(u-1))\,\mathrm{d}u = c^{-\sign
(x)}\alpha^{-1}$, we obtain
\begin{eqnarray*}
\bigl[h(x) \bigr]^{-1}\mathcal{I}_2^{(x)} &=&
\biggl( \frac{1+\beta\sign(x)}{1-\beta\sign(x)} \biggr) \int_{\mathbb
{R}^-} \bigl[ f(xu) -
f(x) \bigr] (-u)^{\alpha-1}\nu \bigl(\sign(x) (u-1) \bigr)\,\mathrm{d}u
\\
& &{} + \frac{ 2\beta\sign(x) }{ 1-\beta\sign(x) } c^{-\sign(x)}\alpha^{-1}f(x).
\end{eqnarray*}
Substituting the values of $a$ and $\beta$ given by (\ref{eqaux1}) and
(\ref{eqaux3}) in the latter equality, it follows
%
\begin{equation}
\label{eq3434} \bigl[h(x) \bigr]^{-1}\mathcal{I}_2^{(x)}
= c^{\sign(x)} \alpha^{-1} I_4^{(x)} -
\alpha^{-1}(\alpha-1)\sign(x)af(x),
\end{equation}
where $I_4^{(x)}$ is the integral
\[
\int_{\mathbb{R}^-} \bigl[ f(xu) - f(x) \bigr] g^\updownarrow(u)
\,\mathrm{d}u.
\]
Thus, the expressions (\ref{eq3232}), (\ref{eq3333}) and (\ref{eq3434}) imply
\begin{eqnarray*}
\mathcal{A}^\updownarrow f(x) &=& |x|^{-\alpha} \bigl[ \bigl(\sign(x)a +
I_2^{(x)} \bigr)xf'(x) + I_1^{(x)}
+ c^{\sign(x)}\alpha^{-1}I_4^{(x)} \bigr]
\\
&&{} + |x|^{-\alpha} \bigl[\alpha^{-1}(\alpha-1)^2
\sign(x)a + I_3^{(x)} \bigr]f(x).
\end{eqnarray*}
Finally, since $h$ is an invariant function for the semigroup of $X^0$,
then $f\equiv1$ belongs to $\mathcal{D}_{\mathcal{A}^\updownarrow}$
and it follows that $\alpha^{-1}(\alpha-1)^2\sign(x)a + I_3^{(x)}=0$.
This ends the proof.
\end{pf*}

\section*{Acknowledgements}

This research has been supported by the ECOS-CONACYT CNRS Research
project M07-M01 and the CNRS-CONACYT International Laboratory Solomon
Lefschetz. Ce travail a b\'en\'efici\'e d'une aide de l'Agence
Nationale de la Recherche portant la r\'ef\'erence ANR-09-BLAN-0084-01.


%


\printhistory


\begin{thebibliography}{16}

\bibitem{Bertoin}
\begin{bbook}[mr]
\bauthor{\bsnm{Bertoin},~\bfnm{Jean}\binits{J.}}
(\byear{1996}).
\btitle{L\'evy Processes}.
\bseries{Cambridge Tracts in Mathematics}
\bvolume{121}.
\baddress{Cambridge}: \bpublisher{Cambridge Univ. Press}.
\bid{mr={1406564}}
\bptok{imsref}%
\end{bbook}
\endbibitem

\bibitem{Bertoin-Yor}
\begin{barticle}[mr]
\bauthor{\bsnm{Bertoin},~\bfnm{Jean}\binits{J.}} \AND
\bauthor{\bsnm{Yor},~\bfnm{Marc}\binits{M.}}
(\byear{2005}).
\btitle{Exponential functionals of {L}\'evy processes}.
\bjournal{Probab. Surv.}
\bvolume{2}
\bpages{191--212}.
\bid{doi={10.1214/154957805100000122}, issn={1549-5787}, mr={2178044}}
\bptok{imsref}%
\end{barticle}
\endbibitem

\bibitem{Caballero-Chaumont}
\begin{barticle}[mr]
\bauthor{\bsnm{Caballero},~\bfnm{M.~E.}\binits{M.E.}} \AND
\bauthor{\bsnm{Chaumont},~\bfnm{L.}\binits{L.}}
(\byear{2006}).
\btitle{Conditioned stable {L}\'evy processes and the {L}amperti
representation}.
\bjournal{J. Appl. Probab.}
\bvolume{43}
\bpages{967--983}.
\bid{doi={10.1239/jap/1165505201}, issn={0021-9002}, mr={2274630}}
\bptok{imsref}%
\end{barticle}
\endbibitem

\bibitem{Caballero-Pardo-Perez10}
\begin{barticle}[mr]
\bauthor{\bsnm{Caballero},~\bfnm{M.~E.}\binits{M.E.}},
\bauthor{\bsnm{Pardo},~\bfnm{J.~C.}\binits{J.C.}} \AND
\bauthor{\bsnm{P{\'e}rez},~\bfnm{J.~L.}\binits{J.L.}}
(\byear{2010}).
\btitle{On {L}amperti stable processes}.
\bjournal{Probab. Math. Statist.}
\bvolume{30}
\bpages{1--28}.
\bid{issn={0208-4147}, mr={2792485}}
\bptok{imsref}%
\end{barticle}
\endbibitem

\bibitem{Caballero-Pardo-Perez11}
\begin{barticle}[mr]
\bauthor{\bsnm{Caballero},~\bfnm{M.~E.}\binits{M.E.}},
\bauthor{\bsnm{Pardo},~\bfnm{J.~C.}\binits{J.C.}} \AND
\bauthor{\bsnm{P{\'e}rez},~\bfnm{J.~L.}\binits{J.L.}}
(\byear{2011}).
\btitle{Explicit identities for {L}\'evy processes associated to symmetric
stable processes}.
\bjournal{Bernoulli}
\bvolume{17}
\bpages{34--59}.
\bid{doi={10.3150/10-BEJ275}, issn={1350-7265}, mr={2797981}}
\bptok{imsref}%
\end{barticle}
\endbibitem

\bibitem{Carmona-Petit-Yor}
\begin{bincollection}[mr]
\bauthor{\bsnm{Carmona},~\bfnm{Philippe}\binits{P.}},
\bauthor{\bsnm{Petit},~\bfnm{Fr{\'e}d{\'e}rique}\binits{F.}} \AND
\bauthor{\bsnm{Yor},~\bfnm{Marc}\binits{M.}}
(\byear{2001}).
\btitle{Exponential functionals of {L}\'evy processes}.
In \bbooktitle{L\'evy Processes}
\bpages{41--55}.
\baddress{Boston, MA}: \bpublisher{Birkh\"auser}.
\bid{mr={1833691}}
\bptok{imsref}%
\end{bincollection}
\endbibitem

\bibitem{Chybiryakov}
\begin{barticle}[mr]
\bauthor{\bsnm{Chybiryakov},~\bfnm{Oleksandr}\binits{O.}}
(\byear{2006}).
\btitle{The {L}amperti correspondence extended to {L}\'evy processes and
semi-stable {M}arkov processes in locally compact groups}.
\bjournal{Stochastic Process. Appl.}
\bvolume{116}
\bpages{857--872}.
\bid{doi={10.1016/j.spa.2005.11.009}, issn={0304-4149}, mr={2218339}}
\bptok{imsref}%
\end{barticle}
\endbibitem

\bibitem{Kiu80}
\begin{barticle}[mr]
\bauthor{\bsnm{Kiu},~\bfnm{Sun~Wah}\binits{S.W.}}
(\byear{1980}).
\btitle{Semistable {M}arkov processes in {$\mathbf{R}^{n}$}}.
\bjournal{Stochastic Process. Appl.}
\bvolume{10}
\bpages{183--191}.
\bid{doi={10.1016/0304-4149(80)90020-4}, issn={0304-4149}, mr={0587423}}
\bptok{imsref}%
\end{barticle}
\endbibitem

\bibitem{Lamperti62}
\begin{barticle}[mr]
\bauthor{\bsnm{Lamperti},~\bfnm{John}\binits{J.}}
(\byear{1962}).
\btitle{Semi-stable stochastic processes}.
\bjournal{Trans. Amer. Math. Soc.}
\bvolume{104}
\bpages{62--78}.
\bid{issn={0002-9947}, mr={0138128}}
\bptok{imsref}%
\end{barticle}
\endbibitem

\bibitem{Lamperti72}
\begin{barticle}[mr]
\bauthor{\bsnm{Lamperti},~\bfnm{John}\binits{J.}}
(\byear{1972}).
\btitle{Semi-stable {M}arkov processes. {I}}.
\bjournal{Z. Wahrsch. Verw. Gebiete}
\bvolume{22}
\bpages{205--225}.
\bid{mr={0307358}}
\bptok{imsref}%
\end{barticle}
\endbibitem

\bibitem{Panti12}
\begin{bmisc}[auto:STB|2012/09/25|13:49:33]
\bauthor{\bsnm{Pant{\'i}},~\bfnm{H.}\binits{H.}}
(\byear{2012}).
\bhowpublished{On L\'evy processes conditioned to avoid zero. Preprint. Available at \arxivurl{arXiv:1304.3191}.}
\bptok{imsref}%
\end{bmisc}
\endbibitem

\bibitem{Revuz-Yor}
\begin{bbook}[mr]
\bauthor{\bsnm{Revuz},~\bfnm{Daniel}\binits{D.}} \AND
\bauthor{\bsnm{Yor},~\bfnm{Marc}\binits{M.}}
(\byear{1999}).
\btitle{Continuous Martingales and {B}rownian Motion},
\bedition{3rd} ed.
\bseries{Grundlehren der Mathematischen Wissenschaften [Fundamental Principles
of Mathematical Sciences]}
\bvolume{293}.
\baddress{Berlin}: \bpublisher{Springer}.
\bid{mr={1725357}}
\bptok{imsref}%
\end{bbook}
\endbibitem

\bibitem{Sato}
\begin{bbook}[mr]
\bauthor{\bsnm{Sato},~\bfnm{Ken-iti}\binits{K.i.}}
(\byear{1999}).
\btitle{L\'evy Processes and Infinitely Divisible Distributions}.
\bseries{Cambridge Studies in Advanced Mathematics}
\bvolume{68}.
\baddress{Cambridge}: \bpublisher{Cambridge Univ. Press}.
\bid{mr={1739520}}
\bptok{imsref}%
\end{bbook}
\endbibitem

\bibitem{Yano10}
\begin{barticle}[mr]
\bauthor{\bsnm{Yano},~\bfnm{Kouji}\binits{K.}}
(\byear{2010}).
\btitle{Excursions away from a regular point for one-dimensional symmetric
{L}\'evy processes without {G}aussian part}.
\bjournal{Potential Anal.}
\bvolume{32}
\bpages{305--341}.
\bid{doi={10.1007/s11118-009-9152-6}, issn={0926-2601}, mr={2603019}}
\bptok{imsref}%
\end{barticle}
\endbibitem

\bibitem{Yor92}
\begin{barticle}[mr]
\bauthor{\bsnm{Yor},~\bfnm{Marc}\binits{M.}}
(\byear{1992}).
\btitle{On some exponential functionals of {B}rownian motion}.
\bjournal{Adv. in Appl. Probab.}
\bvolume{24}
\bpages{509--531}.
\bid{doi={10.2307/1427477}, issn={0001-8678}, mr={1174378}}
\bptok{imsref}%
\end{barticle}
\endbibitem

\end{thebibliography}
\end{document}